\documentclass[review,12pt]{elsarticle}
\usepackage{amsmath}
\usepackage{amsfonts}
\usepackage{geometry}   
\usepackage[parfill]{parskip}    
\usepackage{graphicx}
\usepackage{amssymb}
\usepackage{epstopdf}
\usepackage{psfrag}
\newcommand{\figref}[1]{{Figure~\ref{#1}}}
\newcommand {\bea}{\begin{eqnarray}}
\newcommand {\ea}{\end{eqnarray}}
\newtheorem{proposition}{Proposition}[section]
\newtheorem{theorem}{Theorem}[section]
\newtheorem{Assumption}{Assumption}[section]
\newtheorem{lemma}{Lemma}[section]
\newtheorem{remark}{Remark}[section]

\newenvironment{proof}[1][Proof]{\textbf{#1.} }{\hspace{\stretch{1}}\rule{0.5em}{0.5em}}

\usepackage{xspace}

\usepackage{subfigure}

\newcommand{\thmref}[1]{{Theorem~\ref{#1}}}
\newcommand{\lemref}[1]{{Lemma~\ref{#1}}}

\newcommand{\assref}[1]{{Assumption~\ref{#1}}}
\newcommand{\propref}[1]{{Proposition~\ref{#1}}}

\journal{Computers \& Mathematics with Applications}
\begin{document}
\begin{frontmatter}
\title{A note on exponential Rosenbrock-Euler method for the finite element discretization of a semilinear parabolic partial differential equation}

\author[jdm]{Jean Daniel Mukam}
\ead{jean.d.mukam@aims-senegal.org}
\address[jdm]{Fakult\"{a}t f\"{u}r Mathematik, Technische Universit\"{a}t Chemnitz, 09126 Chemnitz, Germany}

\author[at,atb,atc]{Antoine Tambue}
\cortext[cor1]{Corresponding author}
\ead{antonio@aims.ac.za}
\address[at]{Department of Computing Mathematics and Physics,  Western Norway University of Applied Sciences, Inndalsveien 28, 5063 Bergen.}
\address[atb]{Center for Research in Computational and Applied Mechanics (CERECAM), and Department of Mathematics and Applied Mathematics, University of Cape Town, 7701 Rondebosch, South Africa.}
\address[atc]{The African Institute for Mathematical Sciences(AIMS) of South Africa,
6-8 Melrose Road, Muizenberg 7945, South Africa.}

%

\begin{abstract}

In this paper, we consider the numerical approximation of a general second order semi–linear parabolic
partial differential equation. Equations of this type arise in many contexts, such as transport in porous media.
Using finite element method for space discretization and  the exponential Rosenbrock-Euler method for time discretization, 
we provide a convergence proof in space and time under  only the standard Lipschitz condition of the  nonlinear part, for both smooth and nonsmooth initial solution.
This is in contrast to  restrictive assumptions made  in the literature, where the authors have considered only approximation in time so far in their convergence proofs.
The main result reveals how the convergence orders in both space and time depend heavily on the regularity of the initial data. 
In particular, the method achieves optimal convergence order $\mathcal{O}\left(h^{2}+\Delta t^{2}t_m^{-\eta}\right)$ when the initial data belongs to  the domain of the linear operator.
Numerical simulations to sustain our theoretical result are provided.
\end{abstract}

\end{frontmatter}
\section{Introduction}
\label{intro}
We consider the following abstract Cauchy problem with boundary conditions 
\begin{eqnarray}
\label{model}
\dfrac{du(t)}{dt}= Au(t)+F(u(t)),  \quad u(0)=u_0, \quad t\in(0,T], \quad T>0,
\end{eqnarray}
on the Hilbert  space $H=L^2(\Lambda)$, where $\Lambda$ is an open  subset of $\mathbb{R}^d$ $(d=1,2,3)$, which is supposed to be a convex polygon or has a smooth boundary. 
 The linear operator $A : \mathcal{D}(A)\subset H\longrightarrow H$ is negative, not necessarily self adjoint and generates  an analytic   semigroup $S(t) : = e^{At}$, $t\geq 0$. Without loss of generality, the  nonlinear function $F : H\longrightarrow H$
 is  assumed to be autonomous.  Our main focus  will be on the case where   $A$  is a  general second order  elliptic operator.
 Under some technical conditions (see e.g. \cite{Henry,Pazy}), it is  well  known that the mild solution of \eqref{model} is given by 
 \begin{eqnarray}
 \label{mild1}
 u(t)=S(t)u_0+\int_0^tS(t-s)F\left(u(s)\right)ds, \quad t\in[0,T].
 \end{eqnarray}
  In general, it is hard to find the exact solutions of many PDEs. Numerical approximations are currently  the only  important tools to approximate the solutions.
  Approximations are done at two levels,  spatial approximation and  temporal approximation. The  finite element \cite{Thomee}, finite volume \cite{Antonio3}, finite difference  methods are mostly used for space discretization 
of the problem \eqref{model},  while explicit, semi implicit and fully implicit methods are usually used for time discretization.
References about standard discretization methods for  \eqref{model} can be  found in \cite{Antonio3}.  
Due to the time step size constraints, fully implicit schemes are more popular 
for the time discretization  for quite a long time compared to explicit Euler schemes. However, implicit schemes need at each 
time step a solution of large systems of nonlinear equations. This can be the bottleneck in computations when dealing with realistic problems.
 Recent years, exponential integrators have become an attractive alternative in many evolutions equations \cite{Palencia, Alex1, Alex2, Alex3, Antonio3,Antonio1}.
 Most exponential integrators analyzed  early in  the literature \cite{Palencia,Alex2,Alex3} were bounded on the nonlinear problem as in \eqref{model}
 where the linear part $A$ and the nonlinear function $F$ are explicitly known a priori.
 Such approach is justified in situations where the nonlinear function $F$ is small. Due to the fact that in more realistic applications the nonlinear function $F$ 
 can be stronger\footnote {Typical examples are  semi linear  advection  diffusion reaction equations with stiff reaction term}, Exponential Rosenbrock-Type methods 
 have been proposed in \cite{Alex5,Alex4}, where at every time step, the Jacobian of $F$ is added to the  linear operator $A$. The lower order of them, called 
 Exponential Rosenbrock-Euler method (EREM)  has been proved to be efficient in various applications \cite{Gondal,Antonio1}.
 For smooth  initial solutions, this method is well known to be second order convergence in time \cite{Alex5,Alex4} and have good stability properties in the stochastic context \cite{Antjd1}.
 However in many applications initial solutions are not always smooth. Typical examples are  option pricing in finance  or reaction diffusion advection  with discontinuous initial solution.
 We refer to \cite{Vidar, Thalhammer, Stig, Larsson2,Christian,Alexander2,Alexander1} for standard numerical technique with nonsmooth initial data.
 Recently  exponential Rosenbrock-Euler with nonsmooth initial solution was analysed in \cite{Julia1,Julia2}
 under  the additional hypothesis  \cite[Assumption 1]{Julia1,Julia2}.
  Furthermore, to the best of our knowledge, only convergence in time is investigated for  smooth or nonsmooth initial solution in all existing  Exponential Rosenbrock-Type methods. 
    
  The goal of this paper is to provide a rigorous  convergence proof of EREM in space and time for both smooth and  nonsmooth initial solution under 
  more relaxed conditions than those used in \cite{Julia1,Julia2}.
  Indeed only  the standard  Lipschitz condition of the  nonlinear part  is used in our convergence analysis and optimal convergence orders in space and time  are achieved.
  In fact the method achieves convergence orders of  $\mathcal{O}\left(h^{\beta}+\Delta t^{1+\beta/2}t_m^{-\eta}\right)$, where $\beta$ is the regularity parameter of the initial data (see \assref{Assumption1}) and $\eta$ the parameter defined in \assref{Assumption2}.
  Note that when dealing with space discretization, more novel and careful estimates need to be derived. This is because the constant appearing in the error estimate should not depend on the space  discretization parameter $h$.  
  The space discretization is performed using finite element method. Recent work in \cite{Antonio3} can be used to obtain the similar convergence proof for finite volume method.
 
  The  paper is organized as follows. In  Section \ref{nummethod},  result about the  well posedness are provided along with
  EREM scheme and the main result. The proof of  the main result  is presented in Section \ref{Convergence}. In Section \ref{numericalexperiment}, we present some numerical simulations to sustain our theoretical result.

 \section{Mathematical setting and numerical method}
 \label{nummethod}
 \subsection{Notations, setting and well posedness}
 Let us start by presenting briefly  notations, the main function
spaces and norms that  will be used in this paper. 
We denote by $\Vert \cdot \Vert$ the norm associated to
the inner product $(\cdot ,\cdot )$ of the Hilbert space $H=L^{2}(\Lambda)$. The norms in the Sobolev spaces $H^m(\Lambda),\, m \geqslant 0$ will be denoted by
$\Vert. \Vert_m$. For a Hilbert space $U$ we denote by $\Vert\cdot\Vert_{U}$
the norm of $U$,
$L(U, H)$  the set of bounded linear operators  from
$U$ to $H$.  For ease of notation, we use $L(U,U)=:L(U)$.
In the sequel, for convenience of presentation we take $A$ to be a second-order operator as this simplifies the convergence proof.
More precisely, we assume $A$ to be given by
\begin{eqnarray}
\label{operator}
Au=\sum_{i,j=1}^{d}\dfrac{\partial}{\partial x_i}\left(q_{ij}(x)\dfrac{\partial u}{\partial x_j}\right)-\sum_{i=1}^dq_i(x)\dfrac{\partial u}{\partial x_i},
\end{eqnarray}
where $q_{ij}\in L^{\infty}(\Lambda)$, $q_i\in L^{\infty}(\Lambda)$. We assume that there is a constant $c_1>0$ such that 
\begin{eqnarray}
\label{ellip0}
\sum_{i,j=1}^dq_{ij}(x)\xi_i\xi_j\geq c_1|\xi|^2, \quad  \xi\in \mathbb{R}^d,\quad x\in\overline{\Omega}.
\end{eqnarray}
As in \cite{Suzuki,Antonio2}, we introduce two spaces $\mathbb{H}$ and $V$, such that $\mathbb{H}\subset V$, that depend on the choice of the  boundary conditions for the domain of the operator $A$ and the corresponding bilinear form. For example, for Dirichlet (or first-type) boundary conditions we take 
\begin{eqnarray}
V=\mathbb{H}=H^1_0(\Lambda)=\{v\in H^1(\Lambda) : v=0\quad \text{on}\quad \partial \Lambda\}.
\end{eqnarray}
For Robin (third-type) boundary condition and  Neumann (second-type) boundary condition, which is a special case of Robin boundary condition ($\alpha_0=0$), we take $V=H^1(\Lambda)$
\begin{eqnarray}
\mathbb{H}=\{v\in H^2(\Lambda) : \partial v/\partial v_A+\alpha_0v=0,\quad \text{on}\quad \partial \Lambda\}, \quad \alpha_0\in\mathbb{R}.
\end{eqnarray}
Using  Green's formula and the boundary conditions, we obtain the corresponding bilinear form associated to $-A$,  given by
\begin{eqnarray}
a(u,v)=\int_{\Lambda}\left(\sum_{i,j=1}^dq_{ij}\dfrac{\partial u}{\partial x_i}\dfrac{\partial v}{\partial x_j}+\sum_{i=1}^dq_i\dfrac{\partial u}{\partial x_i}v\right)dx, \quad u,v\in V,
\end{eqnarray}
for Dirichlet and Neumann boundary conditions, and  
\begin{eqnarray}
a(u,v)=\int_{\Lambda}\left(\sum_{i,j=1}^dq_{ij}\dfrac{\partial u}{\partial x_i}\dfrac{\partial v}{\partial x_j}+\sum_{i=1}^dq_i\dfrac{\partial u}{\partial x_i}v\right)dx+\int_{\partial\Lambda}\alpha_0uvdx, \quad u,v\in V.
\end{eqnarray}
for Robin boundary conditions. Using  G\aa rding's inequality, it holds that there exist two positive constants $\lambda_0$ and $c_0$ such that
\begin{eqnarray}
a(v,v)\geq \lambda_0\Vert v \Vert^2_{1}-c_0\Vert v\Vert^2, \quad \forall v\in V.
\end{eqnarray}
By adding and subtracting $c_{0}u $ on the right hand side of (\ref{model}), we obtain a new operator that we still call $A$  corresponding to the new bilinear form that we still call
$a$ such that the following  coercivity property holds
\begin{eqnarray}
\label{ellip}
(-Av,v)=a(v,v)\geq \; \lambda_0\Vert v\Vert_{1}^{2},\;\;\;\;\;\forall v \in V.
\end{eqnarray}
Note that the expression of the nonlinear term $F$ has changed as we included the term $-c_{0}u$
in the new nonlinear term that we still denote by $F$.
The coercivity property (\ref{ellip}) implies that $A$ is sectorial on $L^{2}(\Lambda)$, i.e.  there exist $C_{1}\geq 0$ and $\theta \in (\frac{1}{2}\pi,\pi)$ such that
\begin{eqnarray}
 \Vert (\lambda I -A )^{-1} \Vert_{L(L^{2}(\Lambda))} \leq \dfrac{C_{1}}{\vert \lambda \vert },\;\quad \quad
\lambda \in S_{\theta},
\end{eqnarray}
where $S_{\theta}=\left\lbrace  \lambda \in \mathbb{C} :  \lambda=\rho e^{i \phi},\; \rho>0,\;0\leq \vert \phi\vert \leq \theta \right\rbrace $ (see e.g. \cite{Henry,Larsson2,Antonio2}).
 Therefore  $A$ is the infinitesimal generator of a bounded analytic semigroup $S(t)=e^{t A}$  on $L^{2}(\Lambda)$  such that
\begin{eqnarray}
S(t)= e^{t A}=\dfrac{1}{2 \pi i}\int_{\mathcal{C}} e^{ t\lambda}(\lambda I - A)^{-1}d \lambda,\;\;\;\;\;\;\;
\;t>0,
\end{eqnarray}
where $\mathcal{C}$  denotes a path that surrounds the spectrum of $A $.
The coercivity  property \eqref{ellip} also implies that $-A$ is a positive operator and its fractional powers are well defined  for any $\alpha>0,$ by
\begin{equation}
 \left\{\begin{array}{rcl}
         (-A)^{-\alpha} & =& \frac{1}{\Gamma(\alpha)}\displaystyle\int_0^\infty  t^{\alpha-1}{\rm e}^{tA}dt,\\
         (-A)^{\alpha} & = & ((-A)^{-\alpha})^{-1},
        \end{array}\right.
\end{equation}
where $\Gamma(\alpha)$ is the Gamma function (see \cite{Henry}).

 Throughout this paper, we make the following assumptions, which are  less restrictive  than current  assumptions  used in \cite{Julia1,Julia2}.
 \begin{Assumption}
 \label{Assumption1} 
 The   initial value $u_0\in\mathcal{D}\left((-A)^{\beta/2}\right)$, $\beta\in(0,2]$.
 \end{Assumption}
 \begin{Assumption}
\label{Assumption2}
We assume that the function $F: H\longrightarrow H$  is Lipschitz continuous and twice Fr\'{e}chet differentiable along the strip of the exact solution, i.e. there exists a positive constant $L$ such that 
\begin{eqnarray*}
\Vert F(u)-F(v)\Vert\leq L\Vert u-v\Vert,\quad  u,v\in H,\\
\Vert F_v(v)\Vert_{L(H)}\leq L, \quad \text{and}\quad  \Vert (-A)^{-\eta} F_{vv}(v)\Vert_{ L(H\times H; H)}\leq L, \quad  v\in H,
\end{eqnarray*}
for some $\eta\in\left(\frac{3}{4}, 1\right)$,  where $F_v(v)=D_vF(v):=\dfrac{\partial F}{\partial v}(v)$ and $F_{vv}(v)=D_{vv}F(v):=\dfrac{\partial^2F}{\partial v^2}(v)$.
\end{Assumption}

The following proposition can be found in \cite{Henry}.
\begin{proposition}
\label{proposition1}
Let  $\alpha, \delta\geq 0$ and $0\leq \gamma\leq1$. Then there exists a positive constant $C$ such that the following estimates hold \footnote{\propref{proposition1} also holds   with a uniform constant $C$ (independent of $h$)  when $A$ and $S(t)$ are replaced respectively by their discrete versions $A_h$ and $S_h(t)$ defined in Section \ref{numscheme}, see e.g. \cite{Suzuki,Fredrick,Larsson2,Antonio2,Antjd1,Antonio3,Thomee}.}
\begin{eqnarray*}
\Vert (-A)^{\delta}S(t)\Vert_{L(H)}\leq Ct^{-\delta}, \quad t>0,\quad \Vert (-A)^{\gamma}(\mathbf{I}-S(t))\Vert_{L(H)}\leq Ct^{\gamma},\quad t\geq 0,\\
(-A)^{\delta}S(t)=S(t)(-A)^{\delta}\quad \text{on}\quad \mathcal{D}((-A)^{\delta})\quad \text{and if $\delta\geq \alpha$ then}\quad \mathcal{D}((-A)^{\delta})\subset \mathcal{D}((-A)^{\alpha}).
\end{eqnarray*}
\end{proposition}

The following lemma will be useful in our convergence analysis.
\begin{lemma}
\label{pointulemma}
For any $0\leq\rho\leq 1$ and $0\leq\gamma\leq 2$, there exists a positive constant $C$ such that
\begin{eqnarray}
\label{pointu1}
\int_{t_1}^{t_2}\Vert(-A)^{\rho/2}S(t_2-r)\Vert^2_{L(H)}dr\leq C(t_2-t_1)^{1-\rho},\quad 0\leq t_1\leq t_2\leq T,\\
\label{pointu2}
\left\Vert\int_{t_1}^{t_2} (-A)^{\frac{\gamma}{2}}S(t_2-r)vdr\right\Vert_{L(H)}\leq C(t_2-t_1)^{1-\frac{\gamma}{2}}\Vert v\Vert,\quad 0\leq t_1\leq t_2\leq T,\; v\in H.
\end{eqnarray}
\end{lemma}

\begin{proof} 
The proof of \lemref{pointulemma} for $0\leq \rho<1$ and $0\leq \gamma<2$ is an immediate consequence of \propref{proposition1}. 
The border cases $\rho=1$ and $\gamma=2$ are of special interest in numerical analysis. For instance when analyzing  
an approximation scheme based on finite element method, the convergence order in space depends strongly on the space regularity, 
which is based on \eqref{pointu2}. Therefore a suboptimal space regularity leads to a suboptimal estimate of the convergence order
in space, see e.g. \cite{Thomee} or the discussion in the introduction of \cite{Kruse1}. The proof of \lemref{pointulemma} 
for $\rho=1$ and $\gamma=2$ can also be obtained from \propref{proposition1}, but with a logarithmic loss, which will leads 
to a logarithmic reduction of convergences orders. The proof in the case of self adjoint operator was recently done in \cite[Lemma 3.2]{Kruse1}
and was used in \cite{Kruse2} to achieve optimal convergence order when dealing with stochastic problems. \eqref{pointu1} extends \cite[Lemma 3.2 (iii)]{Kruse1} 
to the case of not necessarily self adjoint operator.   \lemref{pointulemma} allows to achieve optimal convergence order when dealing with not necessarily self adjoint operator.
In fact, let us write $A=A_s+A_n$, where $A_s$ and $A_n$ are respectively the self-adjoint and the non self-adjoint parts of $A$. 
As in \cite[(147)]{Antonio4}, we use the Zassenhaus product formula \cite{Zassen2,Zassen1} to decompose the semigroup $S(t)$ as follows.
\begin{eqnarray}
\label{zas1}
S(t)=e^{At}=e^{(A_s+A_n)t}=e^{A_st}e^{A_nt}\prod_{k=2}^{\infty}e^{C_k},
\end{eqnarray}
where $C_k=C_k(t)$ are called Zassenhaus exponents. In  \eqref{zas1}, let us set 
\begin{eqnarray}
\label{zas2}
S_N(t):=e^{A_nt}\prod_{k=2}^{\infty}e^{C_k},\quad S(t)=S_s(t)S_N(t),
\end{eqnarray}
where $S_s(t):=e^{A_st}$ is the semigroup generated by $A_s$. Using the Baker-Campbell-Hausdorff representation formula \cite{Zassen2,Zassen4,Zassen3}, 
one can prove exactly as in \cite{Antonio4} that $S_N(t)$ is a linear bounded operator. As in \cite{Antonio4} one can prove that $\mathcal{D}\left((-A)^{\alpha})\right)=\mathcal{D}\left((-A_s)^{\alpha}\right)$, $0\leq \alpha\leq 1$. Therefore using \eqref{zas2} and the boundness of $S_N(t)$, it holds that
\begin{eqnarray}
\label{zas3}
\int_{t_1}^{t_2}\left\Vert (-A)^{\rho/2}S(t_2-r)\right\Vert^2_{L(H)}dr&=&\int_{t_1}^{t_2}\left\Vert(-A)^{\rho/2}S_s(t_2-r)S_N(t_2-r)\right\Vert^2_{L(H)}dr\nonumber\\
&\leq& C\int_{t_1}^{t_2}\left\Vert (-A)^{\rho/2}S_s(t_2-r)\right\Vert^2_{L(H)}dr\nonumber\\
&\leq& C\int_{t_1}^{t_2}\left\Vert (-A_s)^{\rho/2}S_s(t_2-r)\right\Vert^2_{L(H)}dr.
\end{eqnarray}
Since $A_s$ is self-adjoint, it follows from \cite[Lemma 3.2 (iii)]{Kruse1} that
\begin{eqnarray*}
\int_{t_1}^{t_2}\left\Vert (-A_s)^{\rho/2}S_s(t_2-r)\right\Vert^2_{L(H)}dr\leq C(t_2-t_1)^{1-\rho}.
\end{eqnarray*}
This completes the proof of \eqref{pointu1}. The proof of \eqref{pointu2} can be found in \cite[Lemma 3.2, (iv)]{Kruse1}, since it is general and does not uses the fact that $A$ is self-adjoint. 
\end{proof}

The well posedness result is given in the following theorem along with optimal regularity results in both space and time.
\begin{theorem}
\label{theorem0}
Under Assumptions \ref{Assumption1},  and  \ref{Assumption2},   the initial value problem \eqref{model} has a unique mild solution $u\in \mathbf{C}([0,T], H)$,  satisfying
\begin{eqnarray}
\label{born_u1}
\Vert u(t)\Vert \leq C(1+\Vert u_0\Vert),\quad \Vert F(u(t))\Vert\leq C(1+\Vert u_0\Vert),\quad t\in[0,T].
\end{eqnarray}
Moreover, the following optimal regularity results in space and time hold
\begin{eqnarray}
\label{spaceregular}
\Vert (-A)^{\beta/2}u(t)\Vert&\leq& C\left(1+\Vert(-A)^{\beta/2}u_0\Vert\right),\quad\quad\quad\quad\quad\quad\quad\quad t\in[0,T],\\
\label{timeregular1}
\Vert u(t_1)-u(t_2)\Vert &\leq& C(t_1-t_2)^{\beta/2}\left(1+\Vert (-A)^{\beta/2}u_0\Vert\right),\quad\quad\quad 0\leq t_1\leq t_2\leq T,
\end{eqnarray}
where $C=C(\beta, T)$ is a positive constant and $\beta$ is the regularity parameter of \assref{Assumption1}.
\end{theorem}

\begin{proof}
For the proof of the existence and the uniqueness, see \cite[Chapter 6, Theorem 1.2, Page 184]{Pazy} or \cite[Theorem 3.29, Page 104]{Gabriel}. 
The proof of  \eqref{born_u1} can be  found in \cite[Theorem 3.29, Page 104]{Gabriel}. Note that the estimates  \eqref{spaceregular} and \eqref{timeregular1} for $\beta=2$ 
 is of  great interest in numerical analysis as they allow to avoid reduction of convergence orders.  Estimates \eqref{spaceregular} and \eqref{timeregular1} can be easily obtained by using the mild form and the regularity estimates of \propref{proposition1}. But this will lead to a reduction of regularity orders for $\beta=2$, which will therefore reduce the convergence orders in time and space when $\beta=2$.  We fill that gap with the help of \lemref{pointulemma}. First of all using \propref{proposition1}, one can easily prove  \eqref{spaceregular} and \eqref{timeregular1} for $\beta\in[0, 2)$.
 
Let us now first prove \eqref{timeregular1}. From \eqref{mild1},  using triangle inequality, it holds that
\begin{eqnarray}
\label{sam3a}
\Vert u(t_2)-u(t_1)\Vert&\leq& \left\Vert \left(e^{At_2}-e^{At_1}\right)u_0\right\Vert+\left\Vert\int_{t_1}^{t_2}e^{A(t_2-s)}F(u(s))ds\right\Vert\nonumber\\
&+&\left\Vert\int_0^{t_1}\left(e^{A(t_2-s)}-e^{A(t_1-s)}\right)F(u(s))ds\right\Vert\nonumber\\
&:=&I_0+I_1+I_2.
\end{eqnarray}
Using \propref{proposition1},  it holds that
\begin{eqnarray}
\label{sam3}
I_0&=&\left\Vert e^{At_1}\left(e^{A(t_2-t_1)}-\mathbf{I}\right)u_0\right\Vert=\left\Vert e^{At_1}\left(e^{A(t_2-t_1)}-\mathbf{I}\right)(-A)^{-1}(-A)u_0\right\Vert\nonumber\\
&\leq& C\left\Vert e^{At_1}\right\Vert_{L(H)}\left\Vert \left(e^{A(t_2-t_1)}-\mathbf{I}\right)(-A)^{-1}\right\Vert_{L(H)}\left\Vert (-A)u_0\right\Vert\nonumber\\
&\leq& C(t_2-t_1)\left\Vert(-A) u_0\right\Vert.
\end{eqnarray}
Using   \propref{proposition1},  it holds that
\begin{eqnarray}
\label{sam4}
I_1\leq\int_{t_1}^{t_2}\left\Vert e^{A(t_2-s)}F(u(s))\right\Vert ds\leq \int_{t_1}^{t_2}\left\Vert e^{A(t_2-s)}\right\Vert_{L(H)}\Vert F(u(s))\Vert ds\leq  C(t_2-t_1).
\end{eqnarray}
Using triangle inequality, we split $I_2$ as follows
\begin{eqnarray}
\label{decompoI2}
I_2&=&\left\Vert\int_0^{t_1}\left(e^{A(t_2-s)}-e^{A(t_1-s)}\right)\left(F(u(s))-F(u(t_1))\right)ds\right\Vert\nonumber\\
&+&\left\Vert \int_0^{t_1}\left(e^{A(t_2-s)}-e^{A(t_1-s)}\right)F(u(t_1))ds\right\Vert=:I_{21}+I_{22}.
\end{eqnarray}
Using \propref{proposition1}, the Lipschitz condition on $F$ and the fact that \eqref{timeregular1} holds for $\beta\in[0,2)$, it follows that 
\begin{eqnarray}
\label{estiI21}
I_{21}&\leq& \int_0^{t_1}\left\Vert \left(e^{A(t_2-t_1)}-\mathbf{I}\right)e^{A(t_1-s)}\left(F(u(s))-F(u(t_1))\right)\right\Vert ds\nonumber\\
&\leq& \int_0^{t_1}\left\Vert \left(e^{A(t_2-t_1)}-\mathbf{I}\right)(-A)^{-1}\right\Vert_{L(H)}\left\Vert(-A) e^{A(t_1-s)}\right\Vert_{L(H)}\nonumber\\
&&\times\left\Vert\left(F(u(s))-F(u(t_1))\right)\right\Vert ds\nonumber\\
&\leq& C(t_2-t_1)\int_0^{t_1}(t_1-s)^{-1}\Vert u(s)-u(t_1)\Vert ds\nonumber\\
&\leq& C(t_2-t_1)\int_0^{t_1}(t_1-s)^{-1+\beta-\epsilon}\left(1+\Vert Au_0\Vert\right)ds\nonumber\\
&\leq& C(t_2-t_1)\left(1+\Vert Au_0\Vert\right). 
\end{eqnarray}
Using   \propref{proposition1} and   \lemref{pointulemma},   it holds  that
\begin{eqnarray}
\label{estiI22}
I_{22}&=& \left\Vert \left(e^{A(t_2-t_1)}-\mathbf{I}\right)\int_0^{t_1}e^{A(t_1-s)}F(u(t_1))ds\right\Vert\nonumber\\
&\leq& \left\Vert \left(e^{A(t_2-t_1)}-\mathbf{I}\right)(-A)^{-1}\right\Vert_{L(H)}\left\Vert\int_0^{t_1}Ae^{A(t_1-s)}F(u(t_1))ds\right\Vert\nonumber\\
&\leq& C(t_2-t_1)\left\Vert\int_0^{t_1}Ae^{A(t_1-s)}F(u(t_1))ds\right\Vert\nonumber\\
&\leq & C(t_2-t_1)\left(1+\Vert Au_0\Vert\right).
\end{eqnarray} 
Substituting \eqref{estiI22} and \eqref{estiI21} in \eqref{decompoI2} yields
\begin{eqnarray}
\label{sam5}
I_2\leq C(t_2-t_1)\left(1+\Vert Au_0\Vert\right).
\end{eqnarray}
Substituting \eqref{sam5}, \eqref{sam4} and \eqref{sam3} in \eqref{sam3a} completes the proof of \eqref{timeregular1}.
 
  Let us now prove \eqref{spaceregular} for $\beta=2$. 
 From the mild form  \eqref{mild1}, it follows by using  triangle inequality,  \lemref{pointulemma} and \propref{proposition1} that
\begin{eqnarray}
\label{es}
\left\Vert -Au(t)\right\Vert&\leq & C\left\Vert -Au_0\right\Vert+\left\Vert\int_0^t -AS(t-s)\left(F(u(s))-F(u(t))\right) ds\right\Vert\nonumber\\
&+&\left\Vert \int_0^t-AS(t-s)F(u(t))ds\right\Vert\nonumber\\
&\leq&  C\left\Vert -Au_0\right\Vert+\int_{0}^t\Vert (-A)S(t-s)\Vert_{L(H)}\Vert F(u(s))-F(u(t))\Vert ds\nonumber\\
&+&C\Vert F(u(t))\Vert\nonumber\\
&\leq& C\left\Vert -Au_0\right\Vert+C\int_0^t(t-s)^{-1}\Vert u(t)-u(s)\Vert ds+C(1+\Vert u(t)\Vert)\nonumber\\
&\leq& C\left(1+\left\Vert -Au_0\right\Vert\right)+C\int_0^t(t-s)^{-\epsilon}ds\leq C\left(1+\left\Vert -Au_0\right\Vert\right).
\end{eqnarray}
 This completes the proof of \eqref{spaceregular}. 
\end{proof}

\subsection{Fully discrete scheme}
\label{numscheme}
For the space  approximation of  problem \eqref{model}, we start by   discretising  our domain $\Lambda$ by a finite triangulation.
Let $\mathcal{T}_h$ be a triangulation with maximal length $h$. Let $V_h \subset V$ denotes the space of continuous and piecewise 
linear functions over the triangulation $\mathcal{T}_h$. We consider the projection $P_h$  defined  from  $H=L^2(\Lambda)$ to $V_h$ by 
\begin{eqnarray}
\label{discrete1}
(P_hu,\chi)=(u,\chi), \quad \forall \chi\in V_h,\, \forall u\in H.
\end{eqnarray}
The discrete operator $A_h : V_h\longrightarrow V_h$ is defined by 
\begin{eqnarray}
\label{discrete2}
(A_h\phi,\chi)=(A\phi,\chi)=-a(\phi,\chi),\quad \forall \phi,\chi\in V_h.
\end{eqnarray}
 As $-A$, the discrete operator $-A_h$  satisfies the coercivity property \eqref{ellip}. Therefore $A_h$  is also a generator 
of a bounded analytic  semigroup $S_h(t) : =e^{tA_h}$, see e.g. \cite{Suzuki,Larsson2,Antonio2}.  As in \cite{Suzuki,Stig,Antonio2}, we characterize the domain of the operator $(-A)^{\beta/2},\, \beta \in \{ 1, 2\}$  as follows.
\begin{eqnarray*}
\mathcal{D}((-A)^{\beta/2})=\mathbb{H}\cap H^{\beta}(\Lambda), \quad \text{ (for Dirichlet boundary conditions)}.\\
\mathcal{D}(-A)=\mathbb{H}, \quad \mathcal{D}((-A)^{1/2})=H^1(\Lambda), \quad \text{(for Robin boundary conditions)},
\end{eqnarray*}
with the following equivlence of norms:
\begin{eqnarray*}
 \quad \Vert v\Vert_{H^r(\Lambda)}\equiv \Vert (-A)^{r/2}v\Vert=:\Vert v\Vert_r,\quad  v\in\mathcal{D}((-A)^{r/2}),\quad r\geq 0.
\end{eqnarray*}
The semi-discrete in space version of problem \eqref{model} consists of  finding $u^h(t)\in V_h$ such that 
\begin{eqnarray}
\label{semi1}
\dfrac{du^h(t)}{dt}=A_hu^h(t)+P_hF(u^h(t)), \quad u^h(0)=P_hu_0,\quad t\in(0,T].
\end{eqnarray}
The operators $A_h$ and $P_hF$ satisfy the same assumptions as $A$ and $F$ respectively. Therefore, Theorem \ref{theorem0} ensures the existence of a unique mild solution  $u^h(t)$ of \eqref{semi1} represented by
\begin{eqnarray}
\label{mild2}
u^h(t)=S_h(t)u^h(0)+\int_0^tS_h(t-s)P_hF(u^h(s))ds,\quad u^h(0)=P_hu_0,\quad t\in(0,T].
\end{eqnarray}
Throughout this paper, without loss of generality, we use  a  fixed time step  $\Delta t=T/M$, $M\in\mathbb{N}$   and we set  $t_m=m\Delta t\in(0,T]$, $m\in \mathbb{N}$. 
For the time discretization, we consider the exponential Rosenbrock-Euler method  to  compute the  numerical approximation $u^h_m$ of $u^h(t_m)$ at discrete time $t_m=m\Delta t \in (0,T]$,  $\Delta t >0$. 
The method is based on the following linearisation of  \eqref{semi1} at each time step
\begin{eqnarray}
\label{semi}
\dfrac{du^h(t)}{dt}=A_hu^h(t)+J_m^hu^h(t)+G^h_m(u^h(t)), \quad t_m\leq t\leq t_{m+1},\quad m=0,\cdots, M-1,
\end{eqnarray}
where $J_m^h$ is the Fr\'{e}chet derivative of $P_hF$ at $u^h_m$ and $G^h_m$ is the remainder  given by
\begin{eqnarray}
\label{remainder}
J^h_m :=D_uP_hF(u^h_m),\quad \quad G^h_m(u^h(t)) :=P_hF(u^h(t))-J_m^hu^h(t).
\end{eqnarray}
Before continuing with the discretization,  let us provide  the following important remarks and lemma.

\begin{remark}
\label{permutation}
Using the properties of the inner product $(.,.)$ and the definition of $P_h$, one can easily prove that $P_h$ is a linear map from $H$ to $V_h$.
Therefore, $D_vP_hv=P_h v$ for all $v\in H$, where $D_v$ is the differential operator (Fr\'{e}chet derivative at $v$). Then it follows that  for all  $ v\in H$ we have 
\begin{eqnarray*}
D_vP_hF(v)=D_v(P_h\circ F)(v)=D_{v}P_h(F(v))\circ D_vF(v)=P_hD_vF(v), \\
D_{vv}(P_hF)(v)=D_v(D_vP_hF(v))=D_v(P_hD_vF(v))=P_hD_{vv}F(v),
\end{eqnarray*}
 where $ f\circ g$  stands for the composition of mappings $f$ and $g$. Therefore
 \begin{eqnarray} 
 \label{swap1}
 J^h_m:=D_uP_hF(u^h_m)=P_hD_uF(u^h_m).
 \end{eqnarray}
 Similarly, for $J^h(v):=D_vP_hF(v)$, the following holds 
 \begin{eqnarray}
 \label{swap2}
 J^h(v)=D_vP_hF(v)=P_hD_vF(v),\quad v\in H.
 \end{eqnarray}
\end{remark}
\begin{remark}
\label{remark1}
Under Assumption \ref{Assumption2}, using  \eqref{swap2} and  the fact that $P_h$ is bounded, 
it follows that the Jacobian satisfies the global Lipschitz condition, i.e. there exists a positive constant $C>0$ such that 
 \begin{eqnarray*}
 \Vert J^h(u)-J^h(v)\Vert_{L(H)}\leq C\Vert u-v\Vert,\quad  u,v\in H.
 \end{eqnarray*}
\end{remark}

\begin{lemma}
\label{lemma2} Under Assumptions  \ref{Assumption1} and \ref{Assumption2},  for all $m\in\mathbb{N}$, 
$A_h+J^h_m$ is a generator of an analytic semigroup $S^h_m(t):=e^{(A_h+J^h_m)t}$, called perturbed semigroup.
 Moreover, $(S^h_m)_{m\in\mathbb{N}}$ is uniformly bounded (independently of $m$ and $h$). 
\end{lemma}
\begin{proof}
 Since $S_h$ is an analytic semigroup, there exist $K\geq 0$  and $w\in\mathbb{R}$ such that 
 \begin{eqnarray*}
 \Vert S_h(t)\Vert_{L(H)}\leq Ke^{wt},  \quad  t\in[0,T].
 \end{eqnarray*}
 Using Assumption \ref{Assumption2} and the fact that $P_h$ is uniformly bounded, it follows by taking the norm in \eqref{swap1} that $J^h_m$ is a  uniformly bounded linear   operator.
 Therefore applying \cite[Chapter 3, Theorem 1.1, page 76]{Pazy} ends the proof. 
\end{proof}

Giving the solution $u^h(t_m)$ at  $t_m$,  applying the variation of constants formula to \eqref{semi} with initial value $u^h(t_m)$ yields  the solution $u^h(t_{m+1})$ at  $t_{m+1}$  in the following mild representation form
\begin{eqnarray}
\label{semi2}
u^h(t_{m+1})=e^{(A_h+J^{h}_m)\Delta t}u^h(t_m)+\int_{t_m}^{t_{m +1}}e^{(A_h+J^{h}_m)(t_{m+1}-s)}G^h_{m}(u^h(s))ds.
\end{eqnarray}
We note that  \eqref{semi2} is the exact solution of \eqref{semi1} at $t_{m+1}$.
To establish  our numerical method, we use the following approximation 
\begin{eqnarray*}
G^h_m(u^h(t_m+s))\approx G^h_m(u^h_m).
\end{eqnarray*}
Therefore the  integral part of \eqref{semi2} can be approximated as follows.
\begin{eqnarray}
\label{constr1}
\int_{t_m}^{t_{m+1}}e^{(A_h+J^h_m)(t_{m+1}-s)}G^h_m(u^h(s))ds&=&\int_{0}^{\Delta t}e^{(A_h+J^h_m)(\Delta t-s)}G^h_m(u^h(t_m+s))ds\nonumber\\
&\approx&(A_h+J^h_m)^{-1}\left(e^{(A_h+J^h_m)\Delta t}-\mathbf{I}\right) G^h_m(u^h_m).
\end{eqnarray}
Inserting \eqref{constr1} in \eqref{semi2} and using the approximation $u^h(t_m)\approx u^h_m$  gives the following approximation $u^h_{m+1}$ of $u^h(t_{m+1})$ at time $t_{m+1}$
\begin{eqnarray}
\label{erem}
u^h_{m+1}= e^{(A_h+J^h_m)\Delta t}u^h_m+(A_h+J^h_m)^{-1}\left(e^{(A_h+J^h_m)\Delta t}-\mathbf{I}\right)G^h_m(u^h_m),\quad m=0, \cdots, M-1.
\end{eqnarray}
The scheme  \eqref{erem}  is called exponential Rosenbrock-Euler method (EREM).
The  numerical scheme \eqref{erem} can be written in the following equivalent form,   efficient   for implementation
\begin{eqnarray*}
u^h_{m+1}=u^h_m+\Delta t\varphi_1\left(\Delta t(A_h+J^h_m)\right)\left[(A_h+J^h_m)u^h_m+G^h_m(u^h_m)\right],
\end{eqnarray*}
where 
\begin{eqnarray*}
\varphi_1(\Delta t(A_h+J^h_m)):=\frac{1}{\Delta t}(A_h+J^h_m)^{-1}\left(e^{(A_h+J^h_m)\Delta t}-\mathbf{I}\right)=\frac{1}{\Delta t}\int_0^{\Delta t}e^{(A_h+J^h_m)(\Delta t-s)}ds.
\end{eqnarray*}
Note that    $\varphi_1\left(\Delta t(A_h+J^h_m)\right)$ is a uniformly bounded operator (see e.g. \cite[Lemma 2.4]{Alex1}).

Having the numerical method in hand, our goal is to examine its  convergence in space and time toward the exact solution in the $L^2(\Lambda)$ norm.
\subsection{Main result}
 Throughout   this paper,  we denote by $C$  any  generic constant independent of $h$, $m$ and $\Delta t$, which may change from one place to another.
 The main result of this paper is formulated in the following theorem.
 \begin{theorem}
\label{mainresult1}
Let  $u$ be  the mild  solution of   problem \eqref{model} and $u^h_m$  its approximation  at time $t_m$ by EREM scheme \eqref{erem}.
Assume that  Assumptions  \ref{Assumption1} and \ref{Assumption2} are fulfilled.   
Then for $m=1,\cdots,M$, it holds that
\begin{eqnarray*}
\Vert u(t_m)-u^h_m\Vert\leq C\left(h^{\beta}+\Delta t^{1+\beta/2}t_m^{-\eta}\right),
\end{eqnarray*}
where $\beta$ is the regularity parameter  from \assref{Assumption1}.
\end{theorem}

\begin{remark}  
Note   that if  the space discretization is performed using finite volume method, recent work in \cite{Antonio3} can be used to obtain similar error estimates with optimal convergence order 1 in space. 
\end{remark}

\section{Proof of the main result}
\label{Convergence}
The proof of the main result need some preparatory results.
\subsection{Preparatory results}
Let us introduce the Ritz representation operator $R_h : V\longrightarrow V_h$ defined by
\begin{eqnarray}
\label{riesz}
(-AR_hv,\chi)=(-Av,\chi)=a(v,\chi),\quad v\in V,\quad \chi\in V_h.
\end{eqnarray}
Under the regularity assumptions on the triangulation and in view of the $V$-ellipticity condition \eqref{ellip0}, it is well known that the following error estimate holds (see e.g. \cite{Suzuki,Larsson2})
\begin{eqnarray}
\label{riesz1}
\Vert R_hv-v\Vert+h\Vert R_hv-v\Vert_{H^1(\Lambda)}\leq Ch^r\Vert v\Vert_{H^r(\Lambda)},\quad v\in V\cap H^r(\Lambda),\quad r\in[1,2].
\end{eqnarray}
Let us consider the following linear problem
\begin{eqnarray}
\label{linearpb}
w'=Aw,\quad t\in(0,T],\quad w(0)=w_0\quad \text{given}. 
\end{eqnarray}
The corresponding semi-discretization in space problem associated to \eqref{linearpb} is: 
\begin{eqnarray}
\label{lineardiscrete}
\text{Find}\quad w_h\in V_h\quad \text{such that} \quad w'_h=A_hw_h,\quad w^0_h=P_hw_0.
\end{eqnarray}
Let us define the following operator 
\begin{eqnarray*}
G_h(t):=S(t)-S_h(t)P_h=e^{-At}-e^{-A_ht}P_h,
\end{eqnarray*}
so that $w(t)-w_h(t)=G_h(t)w_0$. 
The estimate \eqref{riesz1} was used in \cite{Antonio2,Antjd1} to establish the following important lemma, which extends \cite[Theorem 3.5]{Thomee} to the case of not necessary self-adjoint operator $A$.
\begin{lemma}\cite[Lemma 3.1]{Antonio2} 
\label{lemma3}
Let $w$ and $w^h$ be solutions of \eqref{linearpb} and \eqref{lineardiscrete} respectivly.  Assume that $w_0\in\mathcal{D}((-A)^{\alpha/2})$,
then for $r\in[0,2]$ and $0\leq\alpha\leq r$, the following estimates hold
\begin{eqnarray}
\label{deter}
\Vert w(t)-w^h(t)\Vert=\Vert G_h(t)w_0\Vert&\leq& Ch^{r}t^{-(r-\alpha)/2}\Vert w_0\Vert_{\alpha},\quad t\in(0,T].
\end{eqnarray}
\end{lemma}

The following lemma will be useful in our error estimate in space for the nonlinear problem \eqref{model}. It allows to avoid the logarithmic reduction of space order when $\beta=2$.
\begin{lemma}
\label{lemma3a} 
Let \assref{Assumption1} be fulfilled. Let $0\leq\rho\leq1$, then the following estimate holds
\begin{eqnarray}
\label{mil4}
 \left\Vert\int_0^t G_{h}(s)vds\right\Vert\leq Ch^{2-\rho}\Vert v\Vert_{-\rho}, \quad v\in\mathcal{D}((-A)^{-\rho}),\quad t>0.
\end{eqnarray}
\end{lemma}

\begin{proof}
Note that 
\begin{eqnarray*}
 \int_0^tG_h(s)vds&=&\int_0^tA^{-1}AS(s)vds-
 \int_0^tA_h^{-1}A_hS_h(s)P_hvds\nonumber\\
 &=&A^{-1}(S(t)-\mathbf{I})v-A_h^{-1}(S_h(t)-\mathbf{I})P_hv.
\end{eqnarray*}
Therefore 
\begin{eqnarray}
\label{mil5}
\left\Vert\int_0^tG_h(s)vds\right\Vert\leq \Vert(A_h^{-1}P_h-A^{-1})v\Vert+\Vert (S(t)A^{-1}-S_h(t)A_h^{-1}P_h)v\Vert.
\end{eqnarray}
Note that using the definition of $R_h$ and $A_h$ one can easily check that \cite{Antonio2,Larsson2}
\begin{eqnarray}
\label{essai3}
A_hR_h=P_hA.
\end{eqnarray}
Using  \eqref{essai3} and employing \eqref{riesz1} with $r=2-\rho$ yields
\begin{eqnarray}
\label{mil5a}
\Vert (A_h^{-1}P_h-A^{-1})v\Vert&=&\Vert (R_hA^{-1}-A^{-1})v\Vert=\Vert (R_h-\mathbf{I})A^{-1}v\Vert\nonumber\\
&\leq& Ch^{2-\rho}\Vert (-A)^{-1}v\Vert_{2-\rho}= Ch^{2-\rho}\Vert (-A)^{-\rho/2}v\Vert\nonumber\\
&\leq& Ch^{2-\rho}\Vert v\Vert_{-\rho}.
\end{eqnarray}
Using again  \eqref{essai3} and the traingle inequality, it holds that
\begin{eqnarray}
\label{mil6}
\Vert S(t)A^{-1}v-S_h(t)A_h^{-1}P_hv\Vert&=&\Vert S(t)A^{-1}v-S_h(t)R_hA^{-1}v\Vert\nonumber\\
&\leq& \Vert S(t)A^{-1}v-S_h(t)P_hA^{-1}v\Vert\nonumber\\
&+&\Vert S_h(t)(P_hA^{-1}v-R_hA^{-1}v)\Vert.
\end{eqnarray}
Applying \lemref{lemma3} with $r=\alpha=2-\rho$ yields
\begin{eqnarray}
\label{mil7}
\Vert S(t)A^{-1}v-S_h(t)P_hA^{-1}v\Vert\leq Ch^{2-\rho}\Vert (-A)^{-1}v\Vert_{2-\rho}= Ch^{2-\rho}\Vert v\Vert_{-\rho}.
\end{eqnarray}
Using the boundedness of $S_h(t)$,  the triangle inequality,  the best approximation property of the orthogonal projector $P_h$ (see e.g. \cite{Kruse2,Larsson2,Thomee}) and the estimate \eqref{riesz1} with $r=\alpha=2-\rho$, it holds that
\begin{eqnarray}
\label{mil8}
\Vert S_h(t)(P_hA^{-1}v-R_hA^{-1}v)\Vert&\leq& \Vert P_hA^{-1}v-R_hA^{-1}v\Vert\nonumber\\
&\leq& \Vert P_hA^{-1}v-A^{-1}v\Vert+\Vert A^{-1}v-R_hA^{-1}v\Vert\nonumber\\
&=&\Vert (P_h-\mathbf{I})A^{-1}v\Vert+\Vert (R_h-\mathbf{I})A^{-1}v\Vert\nonumber\\
&\leq &\Vert (R_h-\mathbf{I})A^{-1}v\Vert+\Vert (R_h-\mathbf{I})A^{-1}v\Vert\nonumber\\
&=& 2\Vert (R_h-\mathbf{I})A^{-1}v\Vert\nonumber\\
&\leq& Ch^{2-\rho}\Vert (-A)^{-1}v\Vert_{2-\rho}\leq Ch^{2-\rho}\Vert v\Vert_{-\rho}.
\end{eqnarray}
Substituting \eqref{mil8} and \eqref{mil7} in \eqref{mil6} yields
\begin{eqnarray}
\label{mil9}
\Vert S(t)A^{-1}v-S_h(t)A_h^{-1}P_hv\Vert\leq Ch^{2-\rho}\Vert v\Vert_{-\rho}.
\end{eqnarray}
Substituting \eqref{mil9} and \eqref{mil5a} in \eqref{mil5} yields
\begin{eqnarray}
\label{mil10}
\left\Vert \int_0^tG_h(s)vds\right\Vert\leq Ch^{2-\rho}\Vert v\Vert_{-\rho}.
\end{eqnarray}
This completes the proof of the lemma.
\end{proof}

\begin{lemma} 
[Space error]
\label{lemma4}
Let $u(t)$ and $ u^h(t) $ be the mild solutions of \eqref{model} and \eqref{semi1} respectively. 
Assume  that Assumptions  \ref{Assumption1}  and \ref{Assumption2} are fulfilled. Then the following error estimate holds
\begin{eqnarray*}
\Vert u(t)-u^h(t)\Vert \leq Ch^{\beta},
\end{eqnarray*}
where $\beta \in [0,\beta]$ is the regularity parameter  from \assref{Assumption1}.
\end{lemma}

\begin{proof} 
The proof  uses the mild solutions \eqref{mild1} and \eqref{mild2}.  Indeed
\begin{eqnarray}
\label{e}
e(t)&:=&\Vert u(t)-u^h(t)\Vert\nonumber\\
&\leq& \Vert S(t)u_0-S_h(t)P_hu_0\Vert+\left\Vert \int_0^tS(t-s)F(u(s))ds-\int_0^tS_h(t-s)P_hF(u^h(s))ds\right\Vert\nonumber\\
&=:&e_1(t)+e_2(t).
\end{eqnarray}
Using Lemma \ref{lemma3} with   $r=\alpha =\beta$, we obtain 
\begin{eqnarray}
\label{e1}
e_1(t):=\Vert (S(t)-S_h(t)P_h)u_0\Vert\leq Ch^{\beta}\Vert u_0\Vert_{\beta}.
\end{eqnarray}
For the estimation of $e_2(t)$, we use  triangle inequality, the boundedness of $S_h(t-s)$ and \assref{Assumption2} to obtain
\begin{eqnarray}
\label{e21}
e_2(t)&:=& \left\Vert \int_0^tS(t-s)F(u(s))ds-\int_0^tS_h(t-s)P_hF(u^h(s))ds\right\Vert\nonumber\\
&\leq&\int_0^t\Vert S(t-s)F(u(s))-S_h(t-s)P_hF(u^h(s))\Vert ds\nonumber\\
&\leq& \int_0^t\Vert S_h(t-s)P_h(F(u(s))-F(u^h(s)))\Vert ds\nonumber\\
&+&\left\Vert\int_0^t (S(t-s)-S_h(t-s)P_h)F(u(s)) ds\right\Vert\nonumber\\
&\leq& C\int_0^te(s)ds+\left\Vert\int_0^t (S(t-s)-S_h(t-s)P_h)F(u(s)) ds\right\Vert\nonumber\\
&:=& C\int_0^te(s)ds +e_{21}(t).
\end{eqnarray}
To estimate $e_{21}(t)$, we use triangle inequality to obtain
\begin{eqnarray}
\label{re1}
e_{21}(t)&\leq&\left\Vert\int_0^t (S(t-s)-S_h(t-s)P_h)\left(F(u(s))-F(u(t))\right) ds\right\Vert\nonumber\\
&+&\left\Vert\int_0^t (S(t-s)-S_h(t-s)P_h)F(u(t)) ds\right\Vert\nonumber\\
&:=& e_{211}(t)+e_{212}(t).
\end{eqnarray}
Using \lemref{lemma3} with $r=\beta$ and $\alpha=0$, using \assref{Assumption2} and \thmref{theorem0} yields 
\begin{eqnarray}
\label{re2}
e_{211}(t)&\leq& Ch^{\beta}\int_0^t(t-s)^{-\beta/2}\Vert F(u(s))-F(u(t))\Vert ds\nonumber\\
&\leq& Ch^{\beta}\int_0^t(t-s)^{-\beta/2}\Vert u(s)-u(t)\Vert ds\nonumber\\
&\leq & Ch^{\beta}.
\end{eqnarray}
Using \lemref{lemma3a} with $\rho=0$ yields
\begin{eqnarray}
\label{re3}
e_{212}(t)&\leq& Ch^2\Vert F(u(t))\Vert\leq Ch^2\leq Ch^{\beta}.
\end{eqnarray}
Substituting \eqref{re3} and \eqref{re2} in \eqref{re1} yields
\begin{eqnarray}
\label{re4}
e_{21}(t)\leq Ch^{\beta}.
\end{eqnarray}
Substituting \eqref{re4} in \eqref{e21} yields 
\begin{eqnarray}
\label{re5}
e_2(t)\leq Ch^{\beta}+C\int_0^te(s)ds.
\end{eqnarray}
Substituting \eqref{re5} and \eqref{e1} in \eqref{e} yields
\begin{eqnarray}
\label{efin}
e(t)\leq Ch^{\beta}+C\int_0^te(s)ds.
\end{eqnarray}
Applying Gronwall's inequality to \eqref{efin} yields
\begin{eqnarray}
e(t)=\Vert u(t)-u^h(t)\Vert\leq Ch^{\beta}.
\end{eqnarray}
This completes the proof of the lemma.
\end{proof}

\begin{remark}
 \lemref{lemma4} is an improvement of  \cite[Proposition 3.3]{Fredrick} and \cite[Lemma 8]{Antjd1}. In fact for $\beta=2$, 
 there is a logarithmic reduction of order in \cite[Proposition 3.3]{Fredrick} and \cite[Lemma 8]{Antjd1}. This logarithmic 
 reduction also appears in \cite[Theorem 14.3]{Thomee} and \cite[Theorem 1.1]{Stig}. This gap is filled in \lemref{lemma4} 
 with the help of \lemref{lemma3a}. \lemref{lemma3a} can also be used in \cite{Antonio2} to relax the strong regularity assumption
 on $F$ needed to achieve optimal convergence order in space in \cite[Remark 2.9]{Antonio2}.
\end{remark}

\begin{lemma}
\label{lema1}
 Under \assref{Assumption2}, the function $G^h_m$ defined by \eqref{remainder} satisfies the following global Lipschitz  condition
\begin{eqnarray*}
\Vert G^h_m(u^h)-G^h_m(v^h)\Vert\leq C\Vert u^h-v^h\Vert, \quad  m\in\mathbb{N}, \quad  u^h,v^h\in V_h.
\end{eqnarray*}
\end{lemma}
 \begin{proof}
 Using \eqref{remainder}, \assref{Assumption2},  the fact that $P_h$ and $J^h_m$ are uniformly bounded yields
 \begin{eqnarray*}
 \Vert G^h_m(u^h)-G^h_m(v^h)\Vert&\leq&\Vert P_h(F(u^h)-F(v^h))\Vert+\Vert J^h_mu^h-J^h_mv^h\Vert\nonumber\\
 &\leq &C\Vert u^h-v^h\Vert+\Vert J^h_m\Vert_{L(H)}\Vert u^h-v^h\Vert\nonumber\\
 &\leq & C\Vert u^h-v^h\Vert.
 \end{eqnarray*}
 \end{proof}

The proof of the following stability result can be found in \cite[Lemma 4]{Julia1}.
\begin{lemma}
\label{lemma7}
Under  \assref{Assumption1}, the following estimate holds for the perturbed semigroup 
\begin{eqnarray*}
\left\Vert e^{(A_h+J^h_m)\Delta t}\cdots e^{(A_h+J^h_k)\Delta t}\right\Vert_{L(H)}\leq C,\quad  \,\,0\leq k\leq m,
\end{eqnarray*} 
where $C$ is a positive constant independent of $h$, $m$, $k$ and $\Delta t$.

Moreover, for any $\gamma\in[0, 1)$, the following estimate holds
\begin{eqnarray*}
\left\Vert e^{(A_h+J^h_m)\Delta t}\cdots e^{(A_h+J^h_k)\Delta t}(-A_h)^{\gamma}\right\Vert_{L(H)}\leq Ct_{m-k+1}^{-\gamma},\quad  \,\,0\leq k\leq m,
\end{eqnarray*} 
\end{lemma}
\begin{proof}
Let us provide a new proof which does not use any further lemmas,  then simpler than the one in \cite{Julia1}.
Set \begin{eqnarray*}
\left\{\begin{array}{ll}
S^h_{m,k}:=e^{(A_h+J^h_m)\Delta t}e^{(A_h+J^h_{m-1})\Delta t}\cdots e^{(A_h+J^h_k)\Delta t}, \quad \text{if} \quad m\geq k\\
 S^h_{m,k} : =\mathbf{I},\quad \hspace{6cm}\quad \text{if}\quad m<k.
\end{array}
\right.
\end{eqnarray*}
Using a telescopic sum, we expand $S^h_{m,k}$ as follows: 
\begin{eqnarray}
\label{Jul1}
S^h_{m,k}=e^{A_ht_{m+1-k}}+\sum_{j=k}^me^{A_h(t_{m+1}-t_{j+1})}(e^{(A_h+J^h_j)\Delta t}-e^{A_h\Delta t})S^h_{j-1,k}. 
\end{eqnarray}
Taking the norm in both sides of \eqref{Jul1} and using the stability properties of $e^{tA_h}$ yields 
\begin{eqnarray}
\label{Jul2}
\left\Vert S^h_{m,k}\right\Vert_{L(H)}\leq C+C\sum_{j=k}^m\left\Vert e^{(A_h+J^h_j)\Delta t}-e^{A_h\Delta t}\right\Vert_{L(H)}\left\Vert S^h_{j-1,k}\right\Vert_{L(H)}.
\end{eqnarray}
Using the variation of parameter formula (see \cite[(1.2), Page 77]{Pazy}), it holds that 
\begin{eqnarray}
\label{Jul3}
\left(e^{(A_h+J^h_j)\Delta t}-e^{A_h\Delta t}\right)v=\int_0^{\Delta t}e^{A_h(\Delta t-s)}J^h_je^{(A_h+J^h_j)s}vds, \quad  v \in \mathcal{D}(-A).
\end{eqnarray}
Taking the norm in both sides of \eqref{Jul3}, using \propref{proposition1},  \lemref{lemma2} and the fact that  $J^h_j$ is uniformly bounded,  it holds that
\begin{eqnarray}
\label{Jul3a}
\left\Vert \left(e^{(A_h+J^h_j)\Delta t}-e^{A_h\Delta t}\right)v\right\Vert\leq \int_0^{\Delta t}C\Vert v\Vert ds\leq C\Delta t\Vert v\Vert.
\end{eqnarray}
Therefore from \eqref{Jul3a}, we have
\begin{eqnarray}
\label{Jul4}
\left\Vert e^{(A_h+J^h_j)\Delta t}-e^{A_h\Delta t}\right\Vert_{L(H)}\leq C\Delta t.
\end{eqnarray}
Inserting \eqref{Jul4} in \eqref{Jul2} gives
\begin{eqnarray}
\label{Jul5}
\left\Vert S^h_{m,k}\right\Vert_{L(H)}\leq C+C\Delta t\sum_{j=k}^m\left\Vert S^h_{j-1,k}\right\Vert_{L(H)}.
\end{eqnarray}
Applying the discrete Gronwall's lemma to \eqref{Jul5} completes the proof of  Lemma \ref{lemma7}.
\end{proof}

\begin{lemma}
\label{lemma8}
Let  Assumptions \ref{Assumption1} and \ref{Assumption2} be fulfilled. 
\begin{itemize}
\item[(i)] Let $\alpha\in[0,2]$. Then for all $v\in\mathcal{D}((-A)^{\alpha/2}$, it holds that
\begin{eqnarray}
\label{borne1}
\Vert (-A_h)^{\alpha/2}P_hv\Vert\leq C\Vert(-A)^{\alpha/2}v\Vert.
\end{eqnarray}
\item[(ii)] For any $\gamma\in[0,2)$ and $t\in[0,T]$, it holds that
\begin{eqnarray}
\label{guy1}
\Vert u^h(t)\Vert\leq  C\Vert u_0\Vert,\quad \Vert P_hF(u^h(t))\Vert\leq C\Vert u_0\Vert,\\ 
\Vert(-A_h)^{\gamma/2} u^h(t)\Vert\leq C(1+\Vert (-A)^{\gamma/2}u_0\Vert),
\end{eqnarray}
\item[(iii)] For any $v\in H$ and $\alpha\in[0, 1]$, it holds that
\begin{eqnarray*}
\Vert (-A_h)^{-\frac{\eta}{2}}P_hF_{vv}(v)\Vert_{L(H\times H; H)}\leq C.
\end{eqnarray*}
where $u^h(t)$ is the mild solution of \eqref{semi1} represented by \eqref{mild2}.
\end{itemize}
\end{lemma}

\begin{proof}
\begin{itemize}
\item[(i)] The proof of \eqref{borne1} for $0\leq\alpha\leq 1$ with self-adjoint operator can be found in \cite[(2.12)]{Adam}, 
while the case of not necessary self-adjoint operator can be found in \cite[Lemma 1]{Antjd1}. Using  \eqref{discrete2} and the Cauchy-Schwartz inequality, it holds that
\begin{eqnarray}
\label{sam1}
\Vert A_hP_hv\Vert^2=(A_hP_hv,A_hP_hv)=(AP_hv,A_hP_hv)\leq \Vert AP_hv\Vert\Vert A_hP_hv\Vert.
\end{eqnarray} 
It follows from \eqref{sam1} that $\Vert A_hP_hv\Vert\leq \Vert AP_hv\Vert$.
Using the equivalence of norms $\Vert -Aw\Vert\approx \Vert w\Vert_{H^2(\Lambda)}$, $w\in\mathcal{D}(-A)$ (see e.g. \cite{Larsson2}), the fact  that $P_h$ commutes with weak derivatives (see \cite[(28)]{Antjd1}) and the fact $P_h$ is uniformly bounded with respect to $\Vert .\Vert_{L^2(\Lambda)}$, it holds that
\begin{eqnarray}
\label{sam2}
\Vert A_hP_hv\Vert\leq \Vert AP_hv\Vert\leq C\Vert P_hv\Vert_{H^2(\Lambda)}\leq C\Vert v\Vert_{H^2(\Lambda)}\leq C\Vert Av\Vert,\quad v\in \mathcal{D}(-A).
\end{eqnarray}
Inequality \eqref{sam2} shows that \eqref{borne1} holds  for $\alpha=2$. Note that \eqref{borne1} obviously holds for $\alpha=0$. As in \cite{Adam,Kruse2,Larsson2,Antjd1,Antonio3,Thomee}, the intermediate cases follow by the interpolation technique.
\item[(ii)] The proof of (ii) is similar to that of \eqref{born_u1} and \eqref{spaceregular} by using (i).
\item[(iii)] The proof of (iii) can be found in \cite[(70)]{Antonio4}.
\end{itemize}
\end{proof}
\begin{lemma}
\label{lemma9}
Let $u^h(t)$ be the mild solution of \eqref{semi1}. Let Assumptions \ref{Assumption1} and \ref{Assumption2} be fulfilled.
\begin{itemize}
\item[(i)] Then the following estimate holds
\begin{eqnarray}
\label{oster1}
\Vert D_tu^h(t)\Vert\leq Ct^{-1+\beta/2}, \quad  t\in(0,T].
\end{eqnarray}
\item[(ii)] For any $\alpha\in(0,\beta)$, it holds that
\begin{eqnarray}
\label{arrange}
\Vert (-A_h)^{\alpha/2}D_tu^h(t)\Vert\leq Ct^{-1-\alpha/2+\beta/2},\quad t\in(0,T].
\end{eqnarray}
\item[(iii)] The following estimate holds
\begin{eqnarray}
\label{oster3}
 \Vert D^2_tu^h(t)\Vert\leq Ct^{-2+\beta/2}, \quad  t\in(0,T].
\end{eqnarray}
\end{itemize}
where $\beta$ is defined in \assref{Assumption1}.
\end{lemma}

\begin{proof} 
\begin{itemize}
\item[(i)]
  Let us recall that the mild solution $u^h(t)$ satisfies the following semi-discrete problem
\begin{eqnarray}
\label{salut1a}
D_tu^h(t)=A_hu^h(t)+P_hF(u^h(t)),\quad u^h(0)=P_hu_0.
\end{eqnarray}
Therefore $u^h(t)$ is differentiable and its derivative is given by \eqref{salut1a}. Since $A_h$ is a linear operator, it follows that $A_hu^h(t)$ is differentiable. 
The function $P_hF(u^h(t))$ is differentiable as a composition of differentiable maps. Hence $D_tu^h(t)$ is differentiable, i.e. $u^h(t)$ is twice differentiable in time. 
Using the Chain rule and Remark \ref{permutation}, we obtain
\begin{eqnarray}
\label{salut2}
D^2_tu^h(t)=A_hD_tu^h(t)+P_hD_uF(u^h(t))D_tu^h(t).
\end{eqnarray}
 Using the same arguments as above, it follows that $D^3_tu^h(t)$ exists.  As in \cite[Theorem 5.2]{Larsson2}, we set $v^h(t)=tD_tu^h(t)$. Using \eqref{swap2},
 it follows that $v^h(t)$ satisfies the following equation
\begin{eqnarray}
\label{salut3}
D_tv^h(t)=A_hv^h(t)+D_tu^h(t)+P_hD_uF(u^h(t))v^h(t),\quad v^h(0)=0.
\end{eqnarray}
Therefore by Duhamel's principle, we have
\begin{eqnarray}
\label{salut4}
v^h(t)=\int_0^tS_h(t-s)[D_su^h(s)+P_hD_uF(u^h(s))v^h(s)]ds.
\end{eqnarray}
Taking the  norm in both sides of \eqref{salut4}, using the stability properties of $S_h(t-s)$ (see \propref{proposition1}) and the uniformly boundedness of $P_h$  yields
\begin{eqnarray}
\label{salut5}
\Vert v^h(t)\Vert\leq \int_0^t\Vert S_h(t-s)D_su^h(s)\Vert ds+\int_0^t\Vert D_uF(u^h(s))v^h(s)\Vert ds.
\end{eqnarray}
Using \eqref{salut1a}, it holds that
\begin{eqnarray}
\label{salut6}
S_h(t-s)D_su^h(s)=S_h(t-s)A_hu^h(s)+S_h(t-s)P_hF(u^h(s)).
\end{eqnarray}
Taking the norm in both sides of \eqref{salut6}, using \propref{proposition1}, \lemref{lemma8} (i), the boundedness of $P_h$ and \lemref{lemma8} (ii) yields
\begin{eqnarray}
\label{salut7}
\Vert S_h(t-s)D_su^h(s)\Vert&\leq & \Vert S_h(t-s)(-A)^{1-\beta/2}(-A_h)^{\beta/2}u^h(s)\Vert+\Vert P_hF(u^h(s))\Vert\nonumber\\
&\leq& \Vert S_h(t-s)(-A_h)^{1-\beta/2}\Vert_{L(H)}\Vert (-A_h)^{\beta/2}u^h(s)\Vert+C\Vert u_0\Vert\nonumber\\
&\leq & C(t-s)^{-1+\beta/2}\Vert u_0\Vert_{\beta}+C\Vert u_0\Vert\nonumber\\
&\leq & C(t-s)^{-1+\beta/2}.
\end{eqnarray}
Substituting \eqref{salut7} in \eqref{salut5} yields 
\begin{eqnarray}
\label{salut8a}
\Vert v^h(t)\Vert &\leq &C\int_0^t(t-s)^{-1+\beta/2}ds+C\int_0^t\Vert v^h(s)\Vert ds \nonumber\\
&\leq & Ct^{\beta/2}+C\int_0^t\Vert v^h(s)\Vert ds.
\end{eqnarray}
Applying the continuous Gronwall's lemma to \eqref{salut8a} yields
\begin{eqnarray}
\label{salut8}
\Vert v^h(t)\Vert \leq Ct^{\beta/2}.
\end{eqnarray}
Therefore it follows from \eqref{salut8} that
\begin{eqnarray}
\label{salut9}
\Vert D_tu^h(t)\Vert\leq Ct^{-1+\beta/2}.
\end{eqnarray}
This completes the proof (i).
\item[(ii)]
It follows from \eqref{salut4} that
\begin{eqnarray}
\label{salut10}
D_tu^h(t)=t^{-1}\int_0^tS_h(t-s)[D_su^h(s)+sP_hD_uF(u^h(s))D_su^h(s)]ds,\quad t>0.
\end{eqnarray}
Pre-multiplying both sides of \eqref{salut10} by $(-A_h)^{\alpha/2}$ for any $\alpha\in(0,\beta)$ yields 
\begin{eqnarray}
\label{salut11}
&&(-A_h)^{\alpha/2}D_tu^h(t)\nonumber\\
&=&t^{-1}\int_0^t(-A_h)^{\alpha/2}S_h(t-s)[D_su^h(s)+sP_hD_uF(u^h(s))D_su^h(s)]ds.
\end{eqnarray}
Taking the norm in both sides of \eqref{salut11}, using \propref{proposition1}, \assref{Assumption1}, the uniformly boundedness of $P_h$ and \eqref{salut9} yields
\begin{eqnarray}
\label{salut12}
\Vert (-A_h)^{\alpha/2}D_tu^h(t)\Vert& \leq& Ct^{-1}\int_0^t (t-s)^{-1}\left[\Vert D_su^h(s)\Vert+s\Vert P_hD_uF\left(u^h(s)\right)D_su^h(s)\Vert\right] ds\nonumber\\
&\leq& Ct^{-1}\int_0^t(t-s)^{-\alpha/2}[s^{-1+\beta/2}+Cs^{\beta/2}]ds\nonumber\\
&\leq & Ct^{-1}\int_0^t(t-s)^{-\alpha/2}s^{-1+\beta/2}ds\nonumber\\
&\leq & Ct^{-1}t^{-\alpha/2+\beta/2}=Ct^{-1-\alpha/2+\beta/2}.
\end{eqnarray}
\item[(iii)]
We set $w^h(t):=tD_t^2u^h(t)$. Then it holds that
\begin{eqnarray}
\label{bonjour1}
D_tw^h(t)=D^2_tu^h(t)+tD_t^3u^h(t).
\end{eqnarray}
Taking the derivative in both sides of \eqref{salut2}, using the Chain rule and Remark \ref{permutation} yields
\begin{eqnarray}
\label{bonjour2}
D_t^3u^h(t)&=&A_hD^2_tu^h(t)+P_hD_uF(u^h(t))D_t^2u^h(t)\nonumber\\
&+&D_{uu}P_hF(u^h(t))(D_tu^h(t),D_tu^h(t)).
\end{eqnarray}
Substituting \eqref{bonjour2} in \eqref{bonjour1} yields 
\begin{eqnarray}
\label{bonjour3}
D_tw^h(t)&=&A_hw^h(t)+D_tu^h(t)+P_hD_uF(u^h(t))w^h(t)\nonumber\\
&+&tP_hD_{uu}F(u^h(t))(D_tu^h(t),D_tu^h(t)),
\end{eqnarray}
for all $t\in(0,T]$ and $w^h(0)=0$.
Therefore, by Duhamel's principle, it holds that
\begin{eqnarray}
\label{bonjour4}
w^h(t)&=&\int_0^tS_h(t-s)\left[D^2_su^h(s)+P_hD_uF(u^h(s))w^h(s)\right.\nonumber\\
&&+\left. sP_hD_{uu}F(u^h(s))(D_su^h(s),D_su^h(s))\right]ds.
\end{eqnarray}
Taking the norm in both sides of \eqref{bonjour4}, using the uniformly boundedness of $P_h$, \propref{proposition1} and (i) yields
\begin{eqnarray}
\label{bonjour5}
&&\Vert w^h(t)\Vert\nonumber\\
&\leq& \int_0^t\Vert S_h(t-s)D_s^2u^h(s)\Vert ds+C\int_0^t(t-s)^{-\eta}s\Vert D_su^h(s)\Vert^2ds+C\int_0^t\Vert w^h(s)\Vert ds\nonumber\\
&\leq& \int_0^t\Vert S_h(t-s)D_s^2u^h(s)\Vert+C\int_0^t(t-s)^{-\eta}s^{-1+\beta}ds+C\int_0^t\Vert w^h(s)\Vert ds.
\end{eqnarray}
Using \eqref{salut2}, it holds that
\begin{eqnarray}
\label{bonjour6}
S_h(t-s)D_s^2u^h(s)=S_h(t-s)A_hD_su^h(s)+S_h(t-s)P_hD_uF(u^h(s))D_su^h(s).
\end{eqnarray}
Taking the norm in both sides of \eqref{bonjour6},  using \propref{proposition1}, \assref{Assumption2}, \eqref{salut12} and  \eqref{salut9} yields
\begin{eqnarray}
\label{bonjour7}
\Vert S_h(t-s)D^2_su^h(s)\Vert &\leq& \Vert S_h(t-s)A_hD_su^h(s)\Vert +C\Vert D_uF(u^h(s))\Vert_{L(H)}\Vert D_su^h(s)\Vert\nonumber\\
&\leq &\Vert S_h(t-s)A_hD_su^h(s)\Vert+Cs^{-1+\beta/2}\nonumber\\
&=&\Vert D_tS_h(t-s)D_su^h(s)\Vert+Cs^{-1+\beta/2}\nonumber\\
&\leq &C(t-s)^{-1+\alpha/2}\Vert D_su^h(s)\Vert_{\alpha}+Cs^{-1+\beta/2}\nonumber\\
&\leq & C(t-s)^{-1+\alpha/2}s^{-1-\alpha/2+\beta/2}+Cs^{-1+\beta/2}.
\end{eqnarray}
Substituting \eqref{bonjour7} in \eqref{bonjour5} yields
\begin{eqnarray}
\label{bonjour8}
\Vert w^h(t)\Vert&\leq & C\int_0^t(t-s)^{-1+\alpha/2}s^{-1-\alpha/2+\beta/2}ds+C\int_0^t(t-s)^{-\eta}s^{-1+\beta/2}ds\nonumber\\
&+&C\int_0^ts^{-1+\beta}ds+C\int_0^t\Vert w^h(s)\Vert ds.
\end{eqnarray}
Since 
\begin{eqnarray}
\label{bonjour9}
\int_0^t(t-s)^{-1+\alpha/2}s^{-1-\alpha/2+\beta/2}ds\leq Ct^{-1+\beta/2},
\end{eqnarray}
it follows from \eqref{bonjour8} that
\begin{eqnarray}
\label{bonjour10}
\Vert w^h(t)\Vert\leq Ct^{-1+\beta/2}+C\int_0^t\Vert w^h(s)\Vert ds.
\end{eqnarray}
Applying the continuous Gronwall's lemma to \eqref{bonjour10} yields
\begin{eqnarray}
\label{bonjour11}
\Vert w^h(t)\Vert \leq Ct^{-1+\beta/2}.
\end{eqnarray}
It follows therefore from \eqref{bonjour11} that
\begin{eqnarray}
\label{bonjour12}
\Vert D_t^2u^h(t)\Vert \leq Ct^{-2+\beta/2}.
\end{eqnarray} 
\end{itemize}
\end{proof}

\begin{lemma}
\label{lemma10}
Let  Assumptions \ref{Assumption1} and \ref{Assumption2} be fulfilled.
Let   $e^h_m : =u^h(t_m)-u^h_m$ and 
\begin{eqnarray}
\label{bruno}
g_m(t) : = G^h_m(u^h(t))=P_hF(u^h(t))-J^h_mu^h(t).
\end{eqnarray} 
   Then  for all $t_m, t\in(0,T]$, it holds that
\begin{eqnarray*}
\Vert g_m'(t_m)\Vert \leq  Ct_m^{-1+\beta/2} \Vert e^h_m\Vert, \quad \text{and}\quad
\Vert (-A_h)^{-\eta} g_m''(t)\Vert \leq   Ct^{-2+\beta/2},
\end{eqnarray*}
where $\beta$ is defined in \assref{Assumption1} and $\eta$ is defined in \assref{Assumption2}.
\end{lemma}

\begin{proof}
We recall that $J^h_m=D_uP_hF(u^h_m)$ is a linear map. Hence the time derivative  of $J^h_mu^h(t)$ at $t_m$ is given by $J^h_mD_tu^h(t_m)=D_uP_hF(u^h_m)D_tu^h(t_m)$.
Taking the time derivative in \eqref{bruno} and using the Chain rule yields
\begin{eqnarray}
\label{derig1}
g'_m(t)&=&D_uP_hF(u^h(t))D_tu^h(t)-D_uP_hF(u^h_m)D_tu^h(t)\nonumber\\
&=&\left(D_uP_hF(u^h(t))-D_uP_hF(u^h_m)\right)D_tu^h(t).
\end{eqnarray}
 Using  \eqref{swap2}, the fact the projection $P_h$ is bounded and   Remark \ref{remark1}, it follows from \eqref{derig1} that 
\begin{eqnarray*}
\Vert g'_m(t_m)\Vert &\leq& \Vert J^h(u^h(t_m))-J^h(u^h_m)\Vert_{L(H)}\Vert D_tu^h(t_m)\Vert\\
&\leq& C\Vert u^h(t_m)-u^h_m\Vert\Vert D_tu^h(t_m)\Vert= C\Vert e^h_m\Vert\Vert D_tu^h(t_m)\Vert.
\end{eqnarray*}
Using Lemma \ref{lemma9} gives the desired estimate of $\Vert g'_m(t_m)\Vert$. Here the  advantage of the linearisation allows 
to keep $\Vert e^h_m\Vert$ in the upper bound of $\Vert g'_m(t_m)\Vert$ which will be useful in the convergence  proof to reach the  optimal convergence order in time. \\
 Taking the second derivative in \eqref{bruno}, using the chain rule and Remark \ref{permutation} yields 
 \begin{eqnarray}
 \label{derig2}
 g''_m(t)=P_hD_{uu}F(u^h(t))(D_tu^h(t))^2+P_hD_uF(u^h(t))D_t^2u^h(t)-P_hD_uF(u^h_m)D_t^2u^h(t).
 \end{eqnarray}
 Since the projection $P_h$, employing \assref{Assumption2}, it follows from \eqref{derig2} that 
 \begin{eqnarray*}
 \Vert (-A_h)^{-\eta} g''_m(t)\Vert\leq C\Vert D_tu^h(t)\Vert^2+C\Vert D_t^2u^h(t)\Vert.
 \end{eqnarray*}
 Using Lemma \ref{lemma9} completes the proof.
\end{proof}

\subsection{Main proof}
\label{proof_main}
 Let  us  now prove  \thmref{mainresult1}. Using  triangle inequality yields
\begin{eqnarray*}
\Vert u(t_m)-u^h_m\Vert\leq \Vert u(t_m)-u^h(t_m)\Vert+\Vert u^h(t_m)-u^h_m\Vert =: II_1+II_2.
\end{eqnarray*}
 The space error $II_1$ is estimated  by \lemref{lemma4}. It remains to estimate the time error $II_2$. To start, 
we recall that the mild solution  at $t_m$ is given by 
\begin{eqnarray}
\label{new1}
u^h(t_m)&=&e^{(A_h+J^h_{m-1})\Delta t}u^h(t_{m-1})+\int_{t_{m-1}}^{t_m}e^{(A_h+J_{m-1}^h)(t_m-s)}G^h_{m-1}(u^h(s))ds.
\end{eqnarray}
We also recall that the numerical solution \eqref{erem}   at $t_m$ can be written in the following integral form
\begin{eqnarray}
\label{new2}
u^h_m&=&e^{(A_h+J^h_{m-1})\Delta t}u^h_{m-1}+\int_{t_{m-1}}^{t_m}e^{(A_h+J_{m-1}^h)(t_m-s)}G^h_{m-1}(u^h_{m-1})ds.
\end{eqnarray}
If $m=1$, then it follows from \eqref{new1} and \eqref{new2} that 
\begin{eqnarray}
\label{add1a}
II_2:=\Vert u^h(t_1)-u^h_1\Vert=\left\Vert\int_0^{\Delta t}e^{(A_h+J^h_0)(\Delta t-s)}[G^h_0(u^h(s))-G^h_0(u^h_0)]ds \right\Vert.
\end{eqnarray}
Using the uniformly boundedness of $e^{(A_h+J^h_0)t}$ (see \lemref{lemma2}) and \lemref{lema1}, it follows from \eqref{add1a} that 
\begin{eqnarray}
\label{add2a}
II_2&\leq& \int_0^{\Delta t}\left\Vert e^{(A_h+J^h_0)(\Delta t-s)}\right\Vert_{L(H)}\Vert G^h_0(u^h(s))-G^h_0(u^h_0)\Vert ds\nonumber\\
&\leq&C\int_0^{\Delta t}\Vert G^h_0(u^h(s))-G^h_0(u^h_0)\Vert ds=C\int_0^{\Delta t}\Vert G^h_0(u^h(s))-G^h_0(u^h(0))\Vert ds\nonumber\\
&\leq& C\int_0^{\Delta t}\Vert u^h(s)-u^h(0)\Vert ds\nonumber\\
&\leq& C\int_0^{\Delta t}\Vert u^h(s)-u(s)\Vert ds+C\int_0^{\Delta t}\Vert u^h(0)-u(0)\Vert ds+C\int_0^{\Delta t}\Vert u(s)-u(0)\Vert ds\nonumber\\
&\leq& Ch^{\beta}+C\int_0^{\Delta t}s^{\beta/2}ds\leq Ch^{\beta}+C\Delta t^{1+\beta/2}.
\end{eqnarray}
If $m\geq 2$, then iterating the exact solution \eqref{new1}  gives 
\begin{eqnarray}
\label{new1a}
u^h(t_m)&=&e^{(A_h+J^h_{m-1}) \Delta t}e^{(A_h+J^h_{m-2}) \Delta t}\cdots e^{(A_h+J^h_{1}) \Delta t}e^{(A_h+J^h_{0}) \Delta t} u^h(0)\\
&+&\int_{t_{m-1}}^{t_m}e^{(A_h+J^h_{m-1})(t_m-s)}G^h_{m-1}(u^h(s))ds\nonumber\\
&+&\sum_{k=0}^{m-2}\int_{t_{m-k-2}}^{t_{m-k-1}}e^{(A_h+J^h_{m-1}) \Delta t}\cdots e^{(A_h+J^h_{m-k-1}) \Delta t} e^{(A_h+J^h_{m-k-2})(t_{m-k-1}-s)}G^h_{m-k-2}(u^h(s))ds.\nonumber
\end{eqnarray}
For $m\geq 2$,  iterating the numerical solution \eqref{new2}  gives 
\begin{eqnarray}
\label{new2b}
u^h_m&=&e^{(A_h+J^h_{m-1})\Delta t}e^{(A_h+J^h_{m-2}) \Delta t}\cdots e^{(A_h+J^h_{0}) \Delta t} u^h(0)\\
&+&\int_{t_{m-1}}^{t_m}e^{(A_h+J^h_{m-1})(t_{m}-s)}G^h_{m-1}(u^h_{m-1})ds \nonumber\\
&+&\sum_{k=0}^{m-2}\int_{t_{m-k-2}}^{t_{m-k-1}}e^{(A_h+J^h_{m-1}) \Delta t}\cdots e^{(A_h+J^h_{m-k-1}) \Delta t} e^{(A_h+J^h_{m-k-2})(t_{m-k-1}-s)}G^h_{m-k-2}(u^h_{m-k-2})ds\nonumber.
\end{eqnarray}
Therefore, it follows from \eqref{new1a}, \eqref{new2b} and  the triangle inequality that 
\begin{eqnarray*}
\label{new2a}
II_2 &:=&\Vert u^h(t_m)-u^h_m\Vert\\
&\leq&\sum_{k=0}^{m-2} \int_{t_{m-k-2}}^{t_{m-k-1}}\left\Vert e^{(A_h+J^h_{m-1}) \Delta t}\cdots  e^{(A_h+J^h_{m-k-1}) \Delta t}e^{(A_h+J^h_{m-k-2})(t_{m-k-1}-s)}\right.\nonumber\\
&&\left.\left[G^h_{m-k-2}(u^h(s))-G^h_{m-k-2}(u^h_{m-k-2})\right]\right\Vert ds\nonumber\\
&+&\int_{t_{m-1}}^{t_m}\left\Vert e^{(A_h+J^h_{m-1})(t_m-s)}\left[G_{m-1}^h(u^h(s))-G^h_{m-1}(u^h_{m-1})\right]\right\Vert ds\nonumber\\
&\leq& \sum_{k=0}^{m-2}\int_{t_{m-k-2}}^{t_{m-k-1}}\left\Vert e^{(A_h+J^h_{m-1})\Delta t}\cdots e^{(A_h+J^h_{m-k-1})\Delta t}(-A_h)^{\eta}\right\Vert_{L(H)}\nonumber\\
&&\times\Vert (-A_h)^{-\eta}e^{(A_h+J^h_{m-k-2})(t_{m-k-1}-s)}(-A_h)^{\eta}\Vert_{L(H)}\nonumber\\
&&\times \Vert (-A_h)^{-\eta}\left(G^h_{m-k-2}(u^h(s))-G^h_{m-k-2}(u^h_{m-k-2})\right)\Vert ds\nonumber\\
&+&\int_{t_{m-1}}^{t_m}\left\Vert e^{(A_h+J^h_{m-1})(t_{m-1}-s)}\right\Vert_{L(H)}\Vert G^h_{m-1}(u^h(s))-G^h_{m-1}(u^h_{m-1})\Vert ds.
\end{eqnarray*}

Using  Lemmas \ref{lemma7}, \ref{lemma2} and triangle inequality, it holds that
\begin{eqnarray}
\label{new3}
II_2&\leq& C\sum_{k=0}^{m-2}t_{k+1}^{-\eta}\int_{t_{m-k-2}}^{t_{m-k-1}}\Vert (-A_h)^{-\gamma}\left( G^h_{m-k-2}(u^h(s))-G^h_{m-k-2}(u^h_{m-k-2})\right)\Vert ds\nonumber\\
&& +C\int_{t_{m-1}}^{t_m}\Vert G_{m-1}^h(u^h(s))-G^h_{m-1}(u^h_{m-1})\Vert ds\nonumber\\
&\leq& C\sum_{k=0}^{m-2}t_{k+1}^{-\eta}\int_{t_{m-k-2}}^{t_{m-k-1}}\Vert (-A_h)^{-\gamma}\left(G^h_{m-k-2}(u^h(s))-G^h_{m-k-2}(u^h(t_{m-k-2}))\right)\Vert ds\nonumber\\
&+&C\int_{t_{m-1}}^{t_m}\Vert G^h_{m-1}(u^h(s))-G^h_{m-1}(u^h(t_{m-1}))\Vert ds\nonumber\\
&+&C\sum_{k=0}^{m-2}t_{k+1}^{-\eta}\int_{t_{m-k-2}}^{t_{m-k-1}}\Vert (-A_h)^{-\gamma}\left(G^h_{m-k-2}(u^h(t_{m-k-2}))-G^h_{m-k-2}(u^h_{m-k-2})\right)\Vert ds\nonumber\\
&+&C\int_{t_{m-1}}^{t_m}\Vert G^h_{m-1}(u^h(t_{m-1}))-G^h_{m-1}(u^h_{m-1})\Vert ds = : II_{21}+II_{22}+II_{23}+II_{24}.
\end{eqnarray}
Using  \lemref{lema1} and \thmref{theorem0}, it holds that
\begin{eqnarray}
\label{ap}
&&II_{21}+II_{22}\nonumber\\
&=& Ct_m^{-\eta}\int_0^{\Delta t}\Vert G^h_0(u^h(s))-G^h_0(u^h(t_0))\Vert ds\nonumber\\
&+&C\sum_{k=0}^{m-3}t_{k+1}^{-\eta}\int_{t_{m-k-2}}^{t_{m-k-1}}\Vert (-A_h)^{-\gamma}\left(G^h_{m-k-2}(u^h(s))-G^h_{m-k-2}(u^h(t_{m-k-2}))\right)\Vert ds\nonumber\\
&+&C\int_{t_{m-1}}^{t_m}\Vert G^h_{m-1}(u^h(s))-G^h_{m-1}(u^h(t_{m-1}))\Vert ds\nonumber\\
&\leq &Ct_m^{-\eta}\int_0^{\Delta t}\Vert u^h(s)-u^h(0)\Vert ds+C\sum_{k=1}^{m-2}t_{m-k-1}^{-\eta}\int_{t_{k}}^{t_{k+1}}\Vert (-A_h)^{-\gamma}\left(G^h_{k}(u^h(s))-G^h_{k}(u^h(t_{k}))\right)\Vert ds\nonumber\\
&+&C\int_{t_{m-1}}^{t_m}\Vert G^h_{m-1}(u^h(s))-G^h_{m-1}(u^h(t_{m-1}))\Vert ds\nonumber\\
&\leq&Ct_m^{-\eta}\int_0^{\Delta t}\Vert u^h(s)-u(s)\Vert ds+Ct_m^{-\eta}\int_0^{\Delta t}\Vert u^h(0)-u(0)\Vert ds+Ct_m^{-\eta}\int_0^{\Delta t}\Vert u(s)-u(0)\Vert ds\nonumber\\
&+&C\sum_{k=1}^{m-2}t_{m-k-1}^{-\eta}\int_{t_k}^{t_{k+1}}\Vert(-A_h)^{-\gamma}\left(G^h_k(u^h(s))-G^h_k(u^h(t_k))\right)\Vert ds\nonumber\\
&\leq& Ch^{\beta}+C\Delta t^{1+\beta/2}t_m^{-\eta}+C\sum_{k=1}^{m-2}t_{m-k-1}^{-\eta}\int_{t_k}^{t_{k+1}}\Vert(-A_h)^{-\gamma} \left(G^h_k(u^h(s))-G^h_k(u^h(t_k))\right)\Vert ds.
\end{eqnarray}
Using  the fundamental theorem of Analysis and triangle inequality,  we obtain
\begin{eqnarray}
\label{new4}
II_{21} +II_{22}&=&Ch^{\beta}+C\Delta t^{1+\beta/2}t_m^{-\eta}+C\sum_{k=1}^{m-2}t_{m-k-1}^{-\eta}\int_{t_k}^{t_{k+1}}\Vert (-A_h)^{-\gamma}\left(g_k(s)-g_k(t_k)\right)\Vert ds\nonumber\\
&=&Ch^{\beta}+C\Delta t^{1+\beta/2}t_m^{-\eta}+C\sum_{k=0}^{m-2}t_{m-k-1}^{-\eta}\int_{t_k}^{t_{k+1}}\left\Vert\int_{t_k}^s(-A_h)^{-\eta}g'_k(r)dr\right\Vert ds\nonumber\\
&\leq & Ch^{\beta}+C\Delta t^{1+\beta/2}t_m^{-\eta}+ C\sum_{k=1}^{m-2}t_{m-k-1}^{-\eta}\int_{t_k}^{t_{k+1}}\int_{t_k}^s\Vert (-A_h)^{-\eta}g_k'(r)\Vert drds.
\end{eqnarray}
Using again the fundamental theorem of Analysis and triangle inequality yields
\begin{eqnarray}
II_{21}+II_{22}&\leq& Ch^{\beta}+ C\Delta t^{1+\frac{\beta}{2}}t_m^{-\eta}+ C\sum\limits_{k=1}^{m-2}t_{m-k-1}^{-\eta}\int_{t_k}^{t_{k+1}}\int_{t_k}^s\Vert(-A_h)^{-\gamma}\left(g'_k(r)-g'_k(t_k)\right)\Vert drds\nonumber\\
&+&C\sum_{k=1}^{m-2}t_{m-k-1}^{-\eta}\int_{t_k}^{t_{k+1}}\int_{t_k}^s\Vert g'_k(t_k)\Vert drds\nonumber\\
&\leq&Ch^{\beta}+C\Delta t^{1+\beta/2}t_m^{-\eta}+C\sum_{k=1}^{m-2}t_{m-k-1}^{-\eta}\int_{t_k}^{t_{k+1}}\int_{t_k}^s\int_{t_k}^r\Vert (-A_h)^{-\gamma}g''_k(\xi)\Vert d\xi drds\nonumber\\
&+&C\sum_{k=1}^{m-2}t_{m-k-1}^{-\eta}\int_{t_k}^{t_{k+1}}\int_{t_k}^s\Vert g'_k(t_k)\Vert drds.
\end{eqnarray}
Using \lemref{lemma10}, we obtain
\begin{eqnarray}
\label{hoc0c}
II_{21}+II_{22}&\leq& Ch^{\beta}+C\Delta t^{1+\beta/2}t_m^{-\eta}+C\sum_{k=1}^{m-2}t_{m-k-1}^{-\eta}\int_{t_k}^{t_{k+1}}\int_{t_k}^s
\int_{t_k}^r\xi^{-2+\beta/2}d\xi drds\nonumber\\
&+&C\sum_{k=1}^{m-2}t_{m-k-1}^{-\eta}\int_{t_k}^{t_{k+1}}\int_{t_k}^sr^{-1+\beta/2}\Vert e^h_k\Vert drds\nonumber\\
&\leq& Ch^{\beta}+C\Delta t^{1+\beta/2}t_m^{-\eta}+
C\sum_{k=1}^{m-2}t_{m-k-1}^{-\eta}t_k^{-2+\beta/2}\int_{t_k}^{t_{k+1}}\int_{t_k}^s\int_{t_k}^rd\xi drds\nonumber\\
&+&C\sum_{k=1}^{m-1}t_k^{-1+\beta/2}\int_{t_k}^{t_{k+1}}\int_{t_k}^s\Vert e^h_k\Vert drds\nonumber\\
&\leq&Ch^{\beta}+ C\Delta t^{1+\beta/2}t_m^{-\eta}+C\sum_{k=1}^{m-1}t_k^{-2+\beta/2}\Delta t^3+C\sum_{k=1}^{m-2}t_{m-k-1}^{-\eta}t_k^{-1+\beta/2}\Vert e^h_k\Vert\Delta t^2\nonumber\\
&\leq& Ch^{\beta}+C\Delta t^{1+\beta/2}t_m^{-\eta}+C\Delta t^2\sum_{k=1}^{m-2}t_{m-k-1}^{-\eta}t_k^{-2+\beta/2}\Delta t\nonumber\\
&+&C\Delta t\sum_{k=1}^{m-2}t_{m-k-1}^{-\eta}\Vert e^h_k\Vert.
\end{eqnarray}
Let $\lfloor l\rfloor$ be the floor of $l\in \mathbb{N}$. Splitting the sum  in two parts yields
\begin{eqnarray}
\label{cast3b}
&&\Delta t\sum_{k=1}^{m-2}t_{k}^{-2+\frac{\beta}{2}}t_{m-1-k}^{-\eta}
\nonumber\\
&=&\Delta t\sum_{k=1}^{\left\lfloor\frac{m-1}{2}\right\rfloor}t_{k}^{-2+\frac{\beta}{2}}t_{m-k}^{-\eta}+\Delta t\sum_{k=\left\lfloor\frac{m-1}{2}\right\rfloor+1}^{m-1}t_{k}^{-2+\frac{\beta}{2}}t_{m-k}^{-\eta}\nonumber\\
&\leq&\left(\frac{1}{2}t_{m+1}\right)^{-\eta}\Delta t\sum_{k=1}^{m-1}t_{k}^{-2+\frac{\beta}{2}}+\left(\frac{1}{2}t_{m-1}\right)^{-1+\frac{\beta}{2}-\epsilon}\Delta t\sum_{k=1}^{m-1}t_k^{-1+\epsilon}t_{m-k}^{-\eta}\nonumber\\
&\leq& Ct_m^{-\eta}\Delta t\sum_{k=1}^{m-1}t_k^{-2+\frac{\beta}{2}}+C\Delta t^{-1+\frac{\beta}{2}}t_m^{-\epsilon}\Delta t\sum_{k=1}^{m-1}t_k^{-1+\epsilon}t_{m-k}^{-\eta}\nonumber\\
&\leq&  Ct_m^{-\eta}\Delta t\sum_{k=1}^{m-1}t_k^{-2+\frac{\beta}{2}}+C\Delta t^{-1+\frac{\beta}{2}}t_m^{-\eta}.
\end{eqnarray}
Note  that one can easily obtain 
\begin{eqnarray}
\label{cast3bb}
\Delta t\sum_{k=1}^{m-1}t_k^{-2+\frac{\beta}{2}}&=&\Delta t^{-1+\frac{\beta}{2}}\sum_{k=1}^{m-1}k^{-2+\frac{\beta}{2}}.
\end{eqnarray}
The sequence $v_k=k^{-2+\frac{\beta}{2}}$ is decreasing.  Therefore, by comparison with the integral we have
\begin{eqnarray}
\label{cast3cc}
\sum_{k=1}^{m-1}v_k&=&\sum_{k=1}^{m-1}k^{-2+\frac{\beta}{2}}\leq 1+\int_1^mt^{-2+\frac{\beta}{2}}dt\leq C+Cm^{-1+\frac{\beta}{2}}.
\end{eqnarray}
Substituting \eqref{cast3cc} in \eqref{cast3bb} yields
\begin{eqnarray}
\label{cast3c}
\Delta t\sum_{k=1}^{m-1}t_k^{-2+\frac{\beta}{2}}\leq C\Delta t^{-1+\frac{\beta}{2}}+Ct_m^{-1+\frac{\beta}{2}}\leq C\Delta t^{-1+\frac{\beta}{2}}.
\end{eqnarray}
Substituting \eqref{cast3c} in \eqref{cast3b} yields
\begin{eqnarray}
\label{cast3d}
\Delta t\sum_{k=1}^{m-2}t_{k}^{-2+\frac{\beta}{2}}t_{m-k-1}^{-\eta}\leq C\Delta t^{-1+\frac{\beta}{2}}t_m^{-\eta}.
\end{eqnarray}

Substituting \eqref{cast3d} in \eqref{hoc0c} yields
\begin{eqnarray}
\label{hoc4}
II_{21}+II_{22}&\leq& Ch^{\beta}+C\Delta t^{1+\beta/2}t_m^{-\eta}+C\Delta t\sum_{k=1}^{m-2}t_{m-k-1}^{-\eta}\Vert e^h_k\Vert.
\end{eqnarray}
Using Lemma \ref{lema1} we obtain the following estimate for $II_{23}+II_{24}$
\begin{eqnarray}
\label{new5}
II_{23}+II_{24}
&\leq& C\sum_{k=0}^{m-2}\int_{t_{m-k-2}}^{t_{m-k-1}}\Vert u^h(t_{m-k-2})-u^h_{m-k-2}\Vert ds+C\int_{t_{m-1}}^{t_m}\Vert u^h(t_{m-1})-u^h_{m-1}\Vert ds\nonumber\\
&\leq&  C\sum_{k=0}^{m-2}\Delta t \Vert u^h(t_{m-k-2})-u^h_{m-k-2}\Vert+C\Delta t \Vert u^h(t_{m-1})-u^h_{m-1}\Vert \nonumber\\
&\leq & C\Delta t\sum_{k=0}^{m-1}\Vert u^h(t_{k})-u^h_{k}\Vert.
\end{eqnarray}
Inserting \eqref{new5} and \eqref{new4} in \eqref{new3} yields  
\begin{eqnarray}
\label{new6}
II_2=\Vert u^h(t_m)-u^h_m\Vert\leq Ch^{\beta}+C\Delta t^{1+\beta/2}t_m^{-\eta}+C\Delta t\sum_{k=0}^{m-1} t_{m-k-1}^{-\eta}\Vert 
u^h(t_k)-u^h_k\Vert.
\end{eqnarray}
Applying the generalized discrete Gronwall's lemma to \eqref{new6} yields 
\begin{eqnarray}
\label{new7}
II_2=\Vert u^h(t_m)-u^h_m\Vert\leq C\left(h^{\beta}+\Delta t^{1+\beta/2}t_m^{-\eta}\right).
\end{eqnarray}
Combining the estimates of  $II_2$  and \lemref{lemma4} completes the proof of  Theorem \ref{mainresult1}.

\section{Numerical simulations}
\label{numericalexperiment}

Here, we consider  flow and transport in porous media using the  SPE 10 benchmark case data ~\cite{Christie2001a} with the upper 4 layers.
The domain is  $\Lambda=[0,L_{1}]\times [0,L_{2}]\times [0, L_{3}]$.
To deal with high  P\'{e}clet number,  we discretise in space using
the combined finite element-finite volume method, where  the  finite element method is used for diffusion part and the finite volume for advection  part.
The triangulation  $\mathcal{T}$  is  built on a regularity grid  with steps $\Delta x=20$ ft, $\Delta y=10$ ft, and $\Delta z=2$ ft.
The dimensions of  the domain  $\Lambda$ are $L_1=1200$ ft, $L_2=2200$ ft, and $L_3= 8$ ft. The diffusion tensor is $\mathbf{Q}=10^{-4}\mathbf{I}_3=(q_{i,j})$.
We obtain the Darcy velocity field $\mathbf{q}=(q_i)$  by solving the following  system
\begin{equation}
  \label{couple1}
  \nabla \cdot\mathbf{q} =0, \qquad \mathbf{q}=-\mathbf{k} \nabla p.
\end{equation}
%
For pressure  and concentration, we take the Dirichlet boundary condition
\begin{equation*}
\Gamma_{D}=\left\lbrace\left\lbrace 0 \right\rbrace \times \left\lbrace 0 \right\rbrace\times \left[0,L_{3}\right]\right\rbrace \cup\left\lbrace \left\lbrace L_{1} \right\rbrace \times\left\lbrace L_{2}\right \rbrace\times\left[0,L_{3}\right]\right \rbrace,
\end{equation*}
and homogenous Neumann boundary conditions elsewhere  such that
\begin{eqnarray*}
 p&=&\left\lbrace \begin{array}{l}
 3998.96\;\text{ psi } \quad\text{in}\quad\left\lbrace 0 \right\rbrace \times \left\lbrace 0 \right\rbrace\times\left[0,L_{3}\right]\\
\newline\\
 7997.92\;\text{ psi}\quad\text{in}\quad\left\lbrace L_{1} \right\rbrace \times\left\lbrace L_{2} \right\rbrace\times \left[0,L_{3}\right]
 \end{array}\right.
\newline\\
- \mathbf{k} \,\nabla p (x,t)\,\cdot\underbar{n} &=&0\quad\text{in}\quad\Gamma_{N}= \partial \Omega \backslash\Gamma_{D}.
\end{eqnarray*}
Note that in the  SPE 10 benchmark case, the permeability $\mathbf{k}$  diagonal and  highly heterogeneous.
This models a fixed-pressure injector and producer pair located at two diagonally opposite edges of the model, i.e. at $\left\lbrace 0 \right\rbrace \times \left\lbrace 0 \right\rbrace\times\left[0,L_{3}\right]$ and $\left\lbrace L_{1} \right\rbrace \times\left\lbrace L_{2} \right\rbrace\times \left[0,L_{3}\right]$, respectively.

\begin{figure}[h!]
   \begin{center}
  \subfigure[]{
    \label{fig2a}
    \includegraphics[width=0.48\textwidth]{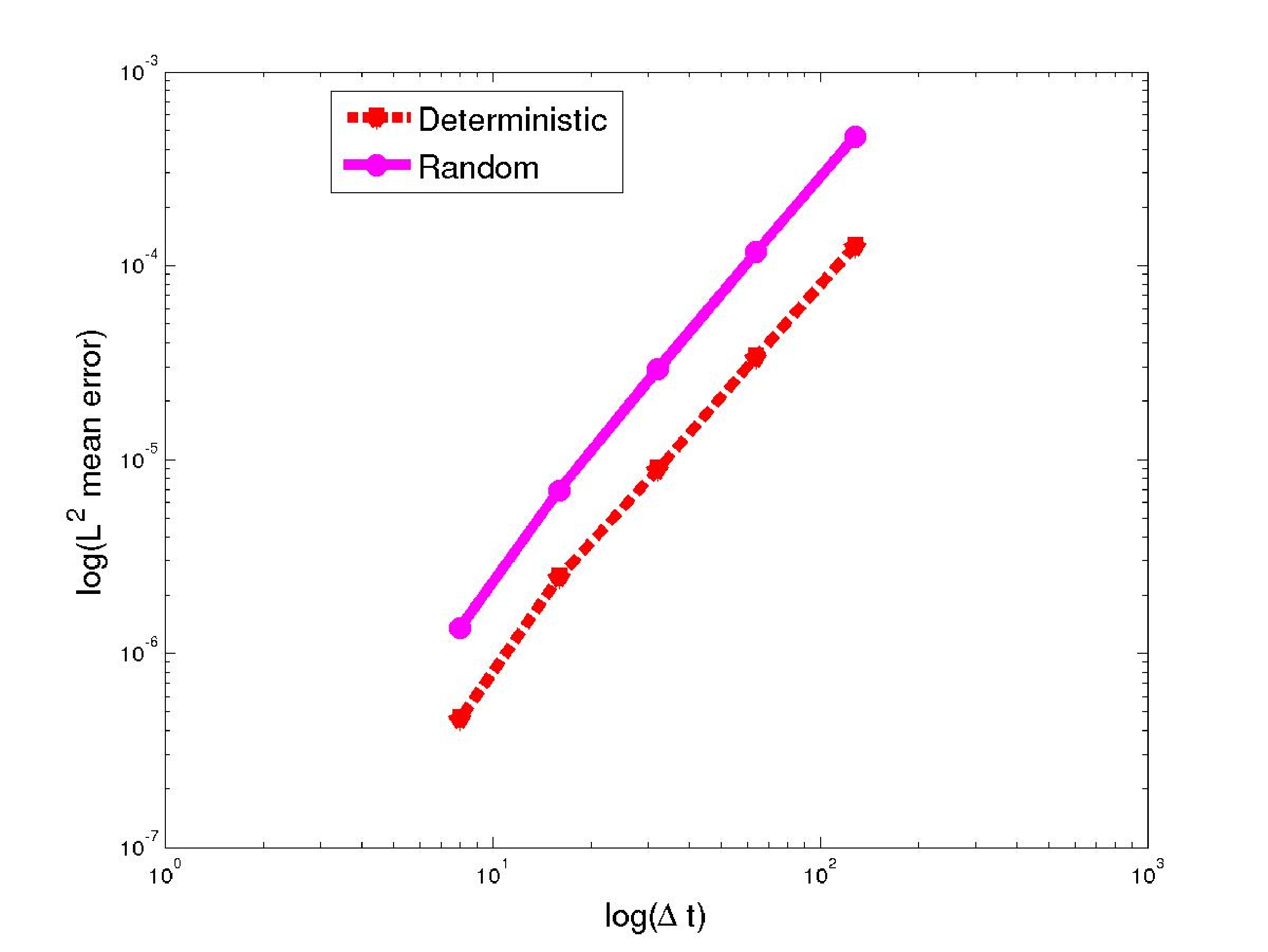}}
  \subfigure[]{
    \label{fig2b}
    \includegraphics[width=1\textwidth]{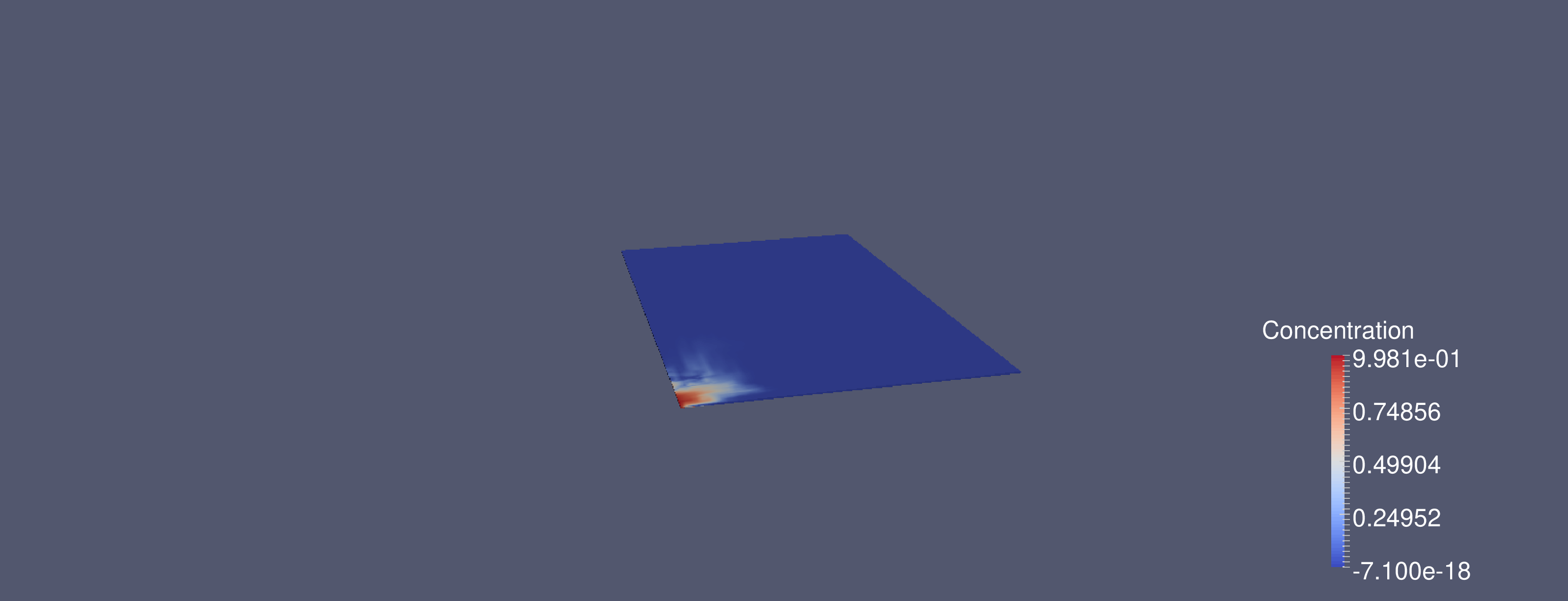}}
   \caption{Graph (a) shows the convergence of the exponential Rosenbrock scheme with deterministic initial  data $u_0=0$ 
   and random initial data  following uniform distribution in the interval $[0,1]$, at the final time $T=8192$ (large time). The convergence orders in time corresponding to the initial value $u_0=0$ and the random initial data are respectively $1.9978$ and $2.0914$, which are in agreement with the theoretical results in \thmref{mainresult1}. The concentration field for the initial solution $u_0=0$ 
   at final  time  $T=8192$ is  presented  in (b).}
   \label{fig2}
\end{center}
\end{figure}

For the concentration, we take
\begin{eqnarray*}
 u&=& 0\quad\text{in}\quad \left\lbrace\left\lbrace 0 \right\rbrace \times \left\lbrace 0 \right\rbrace\times\left[0,L_{3}\right]\right\rbrace\times \left[0,T \right]
 \newline\\
u&=& 1\quad\text{in}\quad \left\lbrace\left\lbrace L_{1} \right\rbrace \times \left\lbrace L_{2} \right\rbrace\times\left[0,L_{3}\right]\right\rbrace\times \left[0,T \right]
\newline\\
-(\mathbf{D}\nabla u)(x,t)\cdot\mathbf{n} &=&0\quad\text{in}\quad\Gamma_{N} \times \left[0,T \right].
\end{eqnarray*}

where $\mathbf{n}$ is the unit outward normal vector to $\Gamma_{N}$. 
For the reaction function we use the classical Langmuir
sorption isotherm  given by $R(u)=(\lambda \beta u)/(1+ \lambda u)$, with $ \lambda =1,\,\beta =10^{-3}$.
In \figref{fig2}, we use the following notations
\begin{itemize}
 \item ''Random'' is used for the numerical solution with random initial data. Indeed  the initial solution  here is not smooth as it 
 follows the uniform distribution in the interval $[0,1]$ and be should be in $L^\infty(\Omega)$. 
 \item ''Deterministic'' is used for numerical solution with null as initial solution.
\end{itemize}
\figref{fig2}(a) shows the convergence of  the exponential Rosenbrock scheme with  both deterministic initial data and random  initial data. 
The orders of convergence are $1.9978$ and $2.0914$ 
respectively.  As the final time  is  large ($T=8192$), the solution is also spatially regular at that time  as we are dealing with  parabolic problem.
Therefore, these convergence  orders are in agreement with our theoretical result in \thmref{mainresult1}.
The final time  is $T=8192$. The concentration field for the numerical solution corresponding  to the  deterministic initial data is presented in \figref{fig2}(b).

\section*{Acknowledgements}
J. D. Mukam was supported by the German Academic Exchange Service (DAAD) (DAAD-Project 57142917) and
A. Tambue was supported by the Robert Bosch Stiftung through the AIMS ARETE CHAIR Programme (Grant No. 11.5.8040.0033.0). 
The authors  thanks Prof. Dr.  Peter Stollmann for his positive and constructive comments. We also thank the referees for their careful reading and useful comments which allowed to improve this paper.


\end{document}